\relax 
\global \edef \sectnotation {{1}}
\global \edef \qsectnotation {{2}}
\global \edef \factsub {{1.1}}
\global \edef \factcolony {{1.3}}
\global \edef \factfewvars {{1.4}}
\global \edef \factdec {{1.5}}
\global \edef \factassocses {{1.6}}
\global \edef \sectsixteen {{2}}
\global \edef \qsectsixteen {{5}}
\global \edef \sectmoreemb {{3}}
\global \edef \qsectmoreemb {{7}}
\global \edef \qcalcc {{8}}
\global \edef \pageminussix {{9}}
\global \edef \thmnox {{3.1}}
\global \edef \sectmoremoreemb {{4}}
\global \edef \qsectmoremoreemb {{12}}
\global \edef \qstep {{14}}
\global \edef \thmthirdemb {{4.1}}
\global \edef \sectreduction {{5}}
\global \edef \qsectreduction {{20}}
\global \edef \thmcolons {{5.2}}
\global \edef \proprad {{5.3}}
\global \edef \thmfirstemb {{5.5}}
\global \edef \qinought {{24}}
\global \edef \thmsecondemb {{5.6}}
\global \edef \thmoneleft {{5.7}}
\global \edef \thmfinal {{5.8}}
\magnification=\magstep1
\baselineskip=15pt
\overfullrule=0pt

\font\footfont=cmr8
\font\footitfont=cmti8
\font\sc=cmcsc10 at 12pt
\font\large=cmbx12 at 14pt

\font\tenmsb=msbm10
\font\sevenmsb=msbm7
\font\fivemsb=msbm5
\newfam\msbfam
\textfont\msbfam=\tenmsb
\scriptfont\msbfam=\sevenmsb
\scriptscriptfont\msbfam=\fivemsb
\def\hexnumber#1{\ifcase#1 0\or1\or2\or3\or4\or5\or6\or7\or8\or9\or
	A\or B\or C\or D\or E\or F\fi}

\mathchardef\subsetneq="2\hexnumber\msbfam28


\def\Ass{\hbox{Ass}\hskip -0.0em}

\def\qed{\hbox{\quad \vrule width 1.6mm height 1.6mm depth 0mm}\vskip 1ex}
\def\eqed{\hbox{\quad \vrule width 1.6mm height 1.6mm depth 0mm}}

\def\D{D}
\def\C{C}
\def\E{E}
\def\F{F}

\def\embprno{$160 n - 270 + 31 d + n(n-1)
+ \delta_{n=2}d^2 + \delta_{n>2}(d^3-d)(n-1)
+ 31(d^{2^1}+ \cdots + d^{2^{n-3}})
+(n-1)d^{2^1} +(n-2) d^{2^2} + \cdots + 3d^{2^{n-3}}
+18 d^{2^{n-2}}$}

\def\mvrule{\vrule height 3.2ex depth 1.6ex width .01em}

\def\elte{c_{02}b_{01}-c_{03}b_{04}}

\newcount\tableco \tableco=0
\def\tabledo{\global\advance\tableco by 1{\the\tableco}}
\def\tabel#1{{\global\edef#1{\the\tableco}}}

\def\label#1{{\global\edef#1{\the\sectno.\the\thmno}}\ignorespaces}
\def\pagelabel#1{\immediate\write\isauxout{\noexpand\global\noexpand\edef\noexpand#1{{\the\pageno}}}{\global\edef#1{\the\pageno}}\ignorespaces}

\def\isnameuse#1{\csname #1\endcsname}

\def\issecond#1#2{#2}
\def\isifundefined#1#2#3{
	\expandafter\ifx\csname #1\endcsname\relax #2
	\else #3 \fi}
\def\pageref#1{\isifundefined{is#1}
	{{\bf ??}\message{Reference `#1' on page [\number\count0] undefined}}
	{\edef\istempa{\isnameuse{is#1}}\expandafter\issecond\istempa\relax}}

\def\today{\ifcase\month\or January\or February\or March\or
  April\or May\or June\or July\or August\or September\or
  October\or November\or December\fi
  \space\number\day, \number\year}

\newcount\sectno \sectno=0
\newcount\thmno \thmno=0
\def \section#1{\vskip 1.2truecm
	\global\advance\sectno by 1 \global\thmno=0
	\noindent{\bf \the\sectno. #1} \vskip 0.6truecm}
\def \thmline#1{\vskip 15pt
	\global\advance\thmno by 1
	\noindent{\bf #1\ \the\sectno.\the\thmno:}\ \ %
	\bgroup \advance\baselineskip by -1pt \it
	\abovedisplayskip =4pt
	\belowdisplayskip =3pt
	\parskip=0pt
	}

\def \thm{\thmline{Theorem}}

\def \endb{\egroup \vskip 1.4ex}
\def \cor{\thmline{Corollary}}

\def \remark{\thmline{Remark}}

\def \fact{\advance\thmno by 1\item{\bf\the\sectno.\the\thmno:}}
\def \prop{\thmline{Proposition}}
\def \proof{\smallskip\noindent {\sl Proof:\ \ }}


\def\label#1{\unskip\immediate\write\isauxout{\noexpand\global\noexpand\edef\noexpand#1{{\the\sectno.\the\thmno}}}
    {\global\edef#1{\the\sectno.\the\thmno}}\unskip\ignorespaces}
\def\sectlabel#1{\immediate\write\isauxout{\noexpand\global\noexpand\edef\noexpand#1{{\the\sectno}}}
    {\global\edef#1{\the\sectno}}}

\def\tabel#1{\hbox to 0pt{\hskip -6em\string#1}\unskip\immediate\write\isauxout{\noexpand\global\noexpand\edef\noexpand#1{{\the\tableco}}}\unskip}

\newwrite\isauxout
\openin1\jobname.aux
\ifeof1\message{No file \jobname.aux}
       \else\closein1\relax\input\jobname.aux
       \fi
\immediate\openout\isauxout=\jobname.aux
\immediate\write\isauxout{\relax}

\ %
\vskip 1ex
\centerline{\large On the embedded primes of the Mayr-Meyer ideals}

\vskip 2ex
\centerline{\sc Irena Swanson}
\footnote{}{The author thanks the NSF for partial support
on grants DMS-0073140 and DMS-9970566.}
\centerline{\sc \today}
\unskip\footnote{ }{{\footitfont 1991 Mathematics Subject Classification.}
13C13, 13P05}
\unskip\footnote{ }{{\footitfont Key words and phrases.}
\footfont Primary decomposition, Mayr-Meyer, complexity of ideals.}

\vskip3ex
\noindent
{\bf Table of contents:}
\newbox\centerdot
\setbox\centerdot=\hbox{\hskip.3em.\hskip.3em}
\def\contentssection#1#2{\vskip.1ex\hbox to 39em{\hskip2em#1\leaders\copy\centerdot\hfil#2}}

\contentssection{Section \sectnotation.\ Notation}{\qsectnotation}
\contentssection{Section \sectsixteen.\ Sixteen embedded components}{\qsectsixteen}
\contentssection{Section \sectmoreemb.\ $15(d+1)$ more embedded primes (plus $d^2$ if $n = 2$)}
{\qsectmoreemb}
\contentssection{Section \sectmoremoreemb.\ $(n-1)(d^3-d)$ more embedded primes, for $n > 2$}
{\qsectmoremoreemb}
\contentssection{Section \sectreduction.\ Reduction to another family}{\qsectreduction}

\vskip 0.5cm

This paper investigates the doubly exponential ideal membership property
of the Mayr-Meyer ideals.
These ideals were first defined by Mayr and Meyer in [MM],
where their doubly exponential behavior was first observed,
and subsequently these ideals were 
further analyzed by Bayer and Stillman [BS], Demazure [D], Koh [K].

The analysis in this paper, as well in as [S1, S2],
is from the point of view of the structure of the associated primes.
The motivation came from a question raised by Bayer, Huneke and Stillman
of whether the doubly exponential behavior is due to the number of minimal
and/or associated primes,
or to the nature of one of them.
The complete answer for the case of the Mayr-Meyer ideals with the fewest possible number
of variables (the case $n = 1$) is given in [S1].
For all other cases,
it was proved in [S2] that the doubly exponential behavior is due to the embedded primes.
[S2] also computed all the minimal components, the minimal primes,
their heights, and the intersection of all the minimal components.

This paper provides partial answers about the embedded primes.
In the analysis a new family of ideals emerges
which also has the doubly exponential ideal membership property.
This new family and its associated primes are further analyzed in [S3].

The main tool used below for finding the associated primes of the Mayr-Meyer
ideals are various short exact sequences,
and the fact that the associated primes of the middle module in a short
exact sequence
is contained in the union of the associated primes of the two other modules.
Theorem~\thmfinal\ gives a set of prime ideals obtained in this way,
which therefore contains all the associated primes of the Mayr-Meyer ideals.
However,
not all the primes in this set need to be associated to the middle module.
Removing the redundant prime ideals is a much harder process,
and is not completed here.
Most of Sections~\sectmoreemb\ and~\sectmoremoreemb\ %
is taken up by removing (some) redundancies.

The Mayr-Meyer ideals $J = J(n,d)$ depend on two parameters,
$n$ and $d$,
where the number of variables in the ring is $O(n)$
and the degree of the given generators of the ideal is $O(d)$.
(See definitions in the first section.)
Both $n$ and $d$ are positive integers.
Throughout this paper it is assumed that $n \ge 2$.

The total number of possibly embedded primes of $J(n,d)$ found in this paper
is \embprno.
This number is doubly exponential in $n$.
Sections~\sectsixteen,  \sectmoreemb, and \sectmoremoreemb\ %
find $31 + 15d + \delta_{n=2}d^2 + \delta_{n>2}(d^3-d)(n-1)$ prime ideals
which are indeed embedded primes of $J(n,d)$,
showing that the number of embedded primes of $J(n,d)$ depends on $n$ and $d$.
Sections~\sectmoreemb\ and \sectmoremoreemb\ %
also prove that there exist no embedded primes of certain kinds.
It is not proved whether $J(n,d)$ in fact has a doubly exponential number
of embedded primes.
On the obtained  list of possibly (but not necessarily) associated primes of $J(n,d)$,
the largest height is achieved by the prime ideals $Q_{23,n-2,n,1,\alpha}$
and $Q_{24}$,
whose heights are 2 less than the dimension of the ring.
However,
I do not know if these ideals are associated.

The generators of the Mayr-Meyer ideals in levels 1 through $n-1$
have similar forms,
so that there is hope that the associated primes of the Mayr-Meyer ideals
could be arrived at via recursion.
I was unable to reduce the search for the associated primes of $J(n,d)$
to that of finding the associated primes of $J(n-1,d)$.
However,
Section~\sectreduction\ modifies the problem
of finding the associated primes of $J(n,d)$
to that of finding the associated primes
of an ideal $K(n,d)$ to which recursion can be applied.
The recursive procedure is carried through in [S3].

Many questions remain about the embedded primes of $J(n,d)$.
Some are listed at the end.

Originally I attempted to find the embedded components,
not just the embedded primes,
but that became unwieldy.
See {\tt http://math.nmsu.edu/\~{}iswanson} for these and other
computations with the Mayr-Meyer ideals which are not included here.

\vskip 3ex

{\bf Acknowledgement.}
I thank Craig Huneke for suggesting this problem
all for all the conversations and enthusiasm for this research.
I also acknowledge the help by computer algebra systems
Macaulay2 and Singular
with which I verified my computations for a few small $n$ and $d$.

\section{Notation}
\sectlabel{\sectnotation}
\pagelabel{\qsectnotation}

The definition below of the Mayr-Meyer ideal is taken from [S2]:
it is somewhat different from the original definition
by Mayr and Meyer in [MM],
but equivalent to the original one from the point of view of primary decompositions.
See [S2] for complete justification.
Namely,
for any fixed integers $n, d \ge 2$,
let $R = k[s, f, b_{ri}, c_{ri}| r = 0, \ldots, n-1; i = 1, \ldots, 4]$
be a polynomial ring in $8n+2$ variables over a field $k$,
and let the Mayr-Meyer ideal $J(n,d)$ be the ideal in $R$
generated by the following polynomials $h_{ri}$:
first the four level 0 generators:
$$
h_{0i} = c_{0i} \left(s -fb_{0i}^d\right), i = 1, 2, 3, 4;
$$
then the eight level 1 generators:
$$
\eqalignno{
& h_{13} = fc_{01} - s c_{02}, \cr
& h_{14} = fc_{04} - s c_{03}, \cr
& h_{15} = s \left(c_{03} - c_{02} \right), \cr
& h_{16} = f \left(c_{02} b_{01} - c_{03} b_{04} \right), \cr
& h_{1,6+i} = f c_{02} c_{1i} \left(b_{02}-b_{1i} b_{03}\right), i = 1, \ldots, 4, \cr
}
$$
the first four level $r$ generators, $r = 2, \ldots, n$:
$$
\eqalignno{
&h_{r3} = s c_{01} c_{11} \cdots c_{r-3,1} \left(
c_{r-2,4} c_{r-1,1} - c_{r-2,1} c_{r-1,2}
\right), \cr
& h_{r4} = s c_{01} c_{11} \cdots c_{r-3,1} \left(
c_{r-2,4} c_{r-1,4} - c_{r-2,1} c_{r-1,3}
\right), \cr
& h_{r5} = s c_{01} c_{11} \cdots c_{r-2,1}
\left( c_{r-1,3}-c_{r-1,2} \right), \cr
& h_{r6} = s c_{01} c_{11} \cdots c_{r-3,1} c_{r-2,4}
\left( c_{r-1,2} b_{r-1,1}-c_{r-1,3} b_{r-1,4} \right), \cr
}
$$
the last four level $r$ generators, $r= 2, ..., n-1$:
$$
h_{r,6+i} = s c_{01} c_{11} \cdots c_{r-3,1} c_{r-2,4} c_{r-1,2} c_{ri}
\left( b_{r-1,2}-b_{ri} b_{r-1,3} \right),
i = 1, \ldots, 4,
$$
and the last level $n$ generator:
$$
h_{n7} = s c_{01} c_{11} \cdots c_{n-3,1} c_{n-2,4} c_{n-1,2}
\left( b_{n-1,2}-b_{n-1,3} \right).
$$

For simpler notation it will be assumed throughout
that the characteristic of $k$ does not divide $d$,
but most of the work goes through without that assumption.
Also, $J(n,d)$ will often be abbreviated to $J$.

For notational purposes we also define the following ideals in $R$:
$$
\eqalignno{
\E &= (s - fb_{01}^d)
+ \left(b_{01} - b_{04}, b_{02}^d - b_{03}^d, b_{01}^d - b_{02}^d\right), \cr
\F &= \left( b_{02} -b_{11}b_{03}, b_{14}-b_{11},b_{13}-b_{11},
b_{12}-b_{11},b_{12}^d - 1\right) \cr
\C_r &= (c_{r1}, c_{r2}, c_{r3}, c_{r4}), r = 0, \ldots, n-1 \cr
\C_n &= (0), \cr
\D_0 &=
\left(c_{04}-c_{01},c_{03}-c_{02},c_{01}-c_{02}b_{01}^d \right), \cr
\D_r &= \left(c_{r4}-c_{r1},c_{r3}-c_{r2},c_{r2}-c_{r1} \right), r = 1, \ldots, n-1, \cr
\D_n &= (0), \cr
B_{0} &= B_{1} = (0), \cr
B_{r} &= \left(1-b_{2i}, 1-b_{3i}, \ldots, 1-b_{ri} |
i = 1, \ldots, 4 \right), r = 2, \ldots, n-1. \cr
B_{kr} &= \left(1-b_{ki}, 1-b_{k+1,i}, \ldots, 1-b_{ri} |
i = 1, \ldots, 4 \right), r = 2, \ldots, n-1. \cr
p_1 &= \C_1 + \E + \D_0, \cr
p_r &= \C_r + \E + \F + B_{r-1} + \D_0 + \D_1 + \cdots + \D_{r-1}, r \ge 2. \cr
}
$$

With this,
here is the table of all the minimal primes over $J(n,d)$,
as computed in [S2],
where $\alpha$ and $\beta$ are $d$th roots of unity,
and $\Lambda$ varies over all the subsets of $\{1,2,3,4\}$:

\vskip 2ex
\halign{\mvrule \hskip 0.1em #\hfil \hskip 0.1em
&& \mvrule \hskip 0.5em \relax #\hfil \mvrule \cr
\noalign{\hrule}
minimal prime & height & component of $J(n,d)$ \cr
\noalign{\hrule}
\noalign{\hrule}
$P_0 = (c_{01}, c_{02}, c_{03}, c_{04})$ & $4$ & $p_0 = P_0$ \cr
\noalign{\hrule}
$P_{1\alpha\beta}= p_1 + (b_{01}- \alpha b_{02}, b_{02}-\beta b_{03})$
& $11$ & $p_{1\alpha\beta}= P_{1\alpha\beta}$ \cr
\noalign{\hrule}
$P_{r\alpha\beta}=p_r$ & $7r + 4$
& $p_{r\alpha\beta}= P_{r\alpha\beta}$, $2 \le r < n$ \cr
\hskip 1em $+ (b_{01}- \alpha b_{02}, b_{02}-\beta b_{03}, \beta-b_{1i})$
& $7n$ & $p_{r\alpha\beta}= P_{r\alpha\beta}$, $r =n$ \cr
\noalign{\hrule}
$P_{-1}= (s,f)$ & $2$ & $p_{-1} = P_{-1}$ \cr
\noalign{\hrule}
$P_{-2}= (s, c_{01}, c_{02}, c_{04}, b_{03}, b_{04})$ & $6$
& \ $p_{-2}= (s, c_{01}, c_{02}, c_{04}, b_{03}^d, b_{04})$ \cr
\noalign{\hrule}
$P_{-3} = (s, c_{01}, c_{04}, b_{02}, b_{03}, c_{02}b_{01}-c_{03}b_{04})$
& $6$ & $p_{-3} = P_{-3}$ \cr
\noalign{\hrule}
$P_{-4\Lambda} =
(s, c_{01}, c_{03}, c_{04}, b_{01}, b_{02})$ & $10$
& $p_{-4\Lambda} =(s, c_{01}, c_{03}, c_{04}, b_{01}, b_{02}^d)$\cr
$\hskip1em + (c_{1i}, b_{1j} |i \not \in \Lambda, j \in \Lambda)$ &
& $+ (c_{1i}|i \not \in \Lambda)$ \cr
& & $+ (b_{1j}^d, b_{02}-b_{1j}b_{03}, b_{1j}-b_{1j'}
|j,j' \in \Lambda)$ \cr
\noalign{\hrule}
}

\vskip 4ex
The intersection of all the components primary to the $P_{-4\Lambda}$
was computed to be
$$
p_{-4}= (s, c_{01}, c_{03}, c_{04}, b_{01}, b_{02}^d)
+ (c_{1i}(b_{02}-b_{1i}b_{03}), c_{1i} b_{1i}^d,
c_{1i}c_{1j}(b_{1i}-b_{1j})| i,j = 1, \ldots, 4).
$$

The following summarizes the elementary facts about primary decompositions
used in the paper:

\vskip2ex\noindent{\bf Facts:}
\bgroup\parindent=3em
\fact
\label{\factsub}
For any ideals $I, I'$ and $I''$ with $I \subseteq I''$,
$(I + I') \cap I'' = I + I' \cap I''$.

\fact
For any ideal $I$ and element $x$,
$(x) \cap I = x (I : x)$.

\fact
\label{\factcolony}
For any ideals $I$ and $I'$,
and any element $x$,
$(I + xI') : x = (I : x) + I'$.

\fact
\label{\factfewvars}
Let $x_1, \ldots, x_n$ be variables over a ring $R$.
Let $S = R[x_1, \ldots, x_n]$.
For any $f_1 \in R$,
$f_2 \in R[x_1]$,
$\ldots$, $f_n \in R[x_1, \ldots, x_{n-1}]$,
let $L$ be the ideal $(x_1 - f_1, \ldots, x_n - f_n)S$ in $S$.
Then an ideal $I$ in $R$ is primary (respectively, prime)
if and only if $IS + L$ is a primary (respectively, prime) in $S$.
Furthermore,
$\cap_i q_i = I$ is a primary decomposition of $I$
if and only if $\cap_i (q_iS + L)$ is a primary decomposition of $IS + L$.

\fact
\label{\factdec}
Let $x$ be an element of a ring $R$ and $I$ an ideal.
Suppose that there is an integer $k$ such that for all $m$,
$I : x^m \subseteq I : x^k$.
Then $I : x^k$ is also denoted as $I : x^\infty$.
Also,
$I = \left(I : x^k\right) \cap \left(I + (x^k)\right)$.
Thus to find a (possibly redundant) primary decomposition of $I$
it suffices to find primary decompositions of (possibly larger) $I : x^k$
and of $I + (x^k)$.

\fact
\label{\factassocses}
Let $I$ be an ideal in a ring $R$.
Then for any $x \in R$,
$\Ass\left({R \over I}\right) \subseteq
\Ass\left({R \over I : x}\right) \cup
\Ass\left({R \over I+ (x)}\right)$,
and every associated prime of
${R \over I : x}$ is an associated prime of ${R \over I}$.
(Use the short exact sequence
$0 \longrightarrow
{R \over I : x} \longrightarrow
{R \over I} \longrightarrow
{R \over I + (x)} \longrightarrow 0$.)

\fact
Let $x_1, \ldots, x_n, y_1, \ldots, y_m$ be variables over a field $k$
and $I$ an ideal in $k[\underline x] = k[x_1, \ldots, x_n]$
and $J$ an ideal in $k[\underline y] = k[y_1, \ldots, y_m]$.
Then
$$
Ik[\underline x,\underline y] \cap Jk[\underline x,\underline y]
= IJk[\underline x,\underline y].
$$
\egroup

We will use the extended Kronecker delta notation $\delta_P$ as follows:
whenever $P$ is true,
then $A \delta_P$ equals $A$,
and when $P$ is false,
$A \delta_P$ has no effect on the rest of the expression.

\section{Sixteen embedded components}
\sectlabel{\sectsixteen}
\pagelabel{\qsectsixteen}

The Mayr-Meyer ideals do have embedded primes.
The (possible) embedded primes will be denoted as $Q_{r\underline{\hbox{\ \ }}}$,
with $r$ varying from $1$ to $24$,
and the second part of the subscript depending on $r$.

Here is the first batch:
for every subset $\Lambda \subseteq \{1, 2,3,4\}$,
define
$$
Q_{1\Lambda}
= \left(s, c_{01},c_{04}, b_{02},b_{03},\elte\right)
+ (c_{1i} | i \not \in \Lambda)
+ (b_{1i}-b_{1j} | i,j \in \Lambda).
$$
Of all the associated primes of $J(n,d)$ found so far,
these primes only contain $P_{-3}$.
We prove below that each of these 16 prime ideals
is an associated prime of $J$,
with its embedded component being
$$
q_{1\Lambda}
= \left(s, c_{01},c_{04}, b_{02}^d,b_{03}^d,\elte\right)
+ (c_{1i} | i \not \in \Lambda)
+ (b_{02}-b_{1i}b_{03},b_{1i}-b_{1j} | i,j \in \Lambda).
$$
It is clear that the sixteen $Q_{1\Lambda}$ are prime
and the sixteen $q_{1\Lambda}$ are primary.
Note that the height of $Q_{1\emptyset}$ is 10,
but if $\Lambda \not = \emptyset$,
then the height of $Q_{1\Lambda}$ equals 9.
Not only is the height of $Q_{1\emptyset}$ larger than that of the others,
but it even contains all $Q_{1\{i\}}$, $i = 1, \ldots, 4$.
By a computation similar to the one of $p_{-4}$ in Section~2 of [S2],
the intersection of all the $q_{1\Lambda}$ is
$$
\eqalignno{
q_1 
&= \left(s, c_{01},c_{04}, b_{02}^d,b_{03}^d,\elte,
c_{1i} (b_{02}-b_{1i}b_{03}),c_{1i}c_{1j}(b_{1i}-b_{1j}) | i,j =1, \ldots, 4
\right). \cr
}
$$
Observe that $J : (fc_{02}c_{03}(c_{02}-c_{03}))^\infty$ contains
$$
\left(s, c_{01},c_{04}, b_{02}^d,b_{03}^d,\elte,
c_{1i} (b_{02}-b_{1i}b_{03}) \right).
$$
This latter ideal contains $J$ and decomposes as
$$
\eqalignno{
&=
\bigcap_{\Lambda} \left(s, c_{01},c_{04}, b_{02}^d,
b_{03}^d,\elte,
c_{1i}, b_{02}-b_{1j}b_{03} | i \not \in \Lambda, j \in \Lambda \right) \cr
&=
\left(s, c_{01},c_{04},b_{02},b_{03},\elte\right) \cr
&\hskip1em \cap
\bigcap_{\Lambda} \left(s, c_{01},c_{04},b_{02}^d,
b_{03}^d,\elte,c_{1i}, b_{02}-b_{1j}b_{03},b_{1j}-b_{1j'}
| i \not \in \Lambda, j,j' \in \Lambda \right) \cr
&= p_{-3} \cap q_1. \cr
}
$$
The element $fc_{02}c_{03}(c_{02}-c_{03})$ is a non-zerodivisor on these components,
which proves that
$$
J : (fc_{02}c_{03}(c_{02}-c_{03}))^\infty= p_{-3} \cap q_1.
$$
To prove that each $Q_{1\Lambda}$ is associated to $J$,
it now suffices to prove that none of the $Q_{1\Lambda}$-primary components
of $J : (fc_{02}c_{03}(c_{02}-c_{03}))^\infty$ is redundant.
So let $\Lambda$ be a subset of $\{1,2,3,4\}$.
First suppose that $\Lambda \not = \emptyset$.
Then
$J' = J : (fc_{02}c_{03}(c_{02}-c_{03}))^\infty
\bigl(\prod_{{i \in \Lambda} \atop {j \not \in \Lambda}}
c_{1i}(b_{1i}-b_{1j})\bigr)$
is exactly $p_{-3} \cap q_{1\Lambda} \not = p_{-3}$,
so that $Q_{1\Lambda}$ is associated to $J$.

Finally suppose that $\Lambda = \emptyset$.
Then
$J' = J : (fc_{02}c_{03}(c_{02}-c_{03}))^\infty
\bigl(\prod_{i\not = j}(b_{1i}-b_{1j})\bigr)$
is exactly $p_{-3} \cap q_{1\emptyset} \cap \bigcap_{i=1}^4 q_{1\{i\}}$.
The element $b_{02}^{d-1}b_{03} \prod_i b_{1i}$
is in $p_{-3} \cap \bigcap_{i=1}^4 q_{1\{i\}}$
but not in $q_{1\emptyset}$,
which proves that $Q_{1\emptyset}$ is associated to $J$.

Thus $J(n,d)$ has at least 16 embedded primes,
which are as follows:

\vskip 2ex
\halign{\mvrule \hskip 0.1em #\hfil \hskip 0.1em
&& \mvrule \hskip 0.5em \relax #\hfil \mvrule \cr
\noalign{\hrule}
embedded prime & height & component of $J(n,d)$ \cr
\noalign{\hrule}
\noalign{\hrule}
$Q_{1\Lambda}
= \left(s, c_{01},c_{04}, b_{02},b_{03}\right)$
& 9, if $\Lambda \not = \emptyset$
& $q_{1\Lambda}
= \left(s, c_{01},c_{04}, b_{02}^d,b_{03}^d\right)$ \cr
\hskip1em $+ (\elte) + (c_{1i} | i \not \in \Lambda)$
& 10, if $\Lambda = \emptyset$
& \hskip1em$+ (\elte) + (c_{1i} | i \not \in \Lambda)$ \cr
\hskip1em $+ (b_{1i}-b_{1j} | i,j \in \Lambda)$
& & \hskip1em $+ (b_{02}-b_{1i}b_{03},b_{1i}-b_{1j} | i,j \in \Lambda)$ \cr
\noalign{\hrule}
}

\section{
$\bf 15(d+1)$ more embedded primes (plus $\bf d^2$ if $\bf n = 2$)}
\sectlabel{\sectmoreemb}
\pagelabel{\qsectmoreemb}

In this section we find $15(d+1) + d^2$ more embedded primes of $J(n,d)$.
This shows that the number of embedded primes of $J(n,d)$ depends on $d$.
As usual,
$\Lambda$ ranges over subset of $\{1,2,3,4\}$
and $\alpha, \beta$ over the $d$th roots of unity:
$$
\eqalignno{
Q_{2\Lambda\alpha} &=
\left(s, c_{01}, c_{03}-c_{02}, c_{04}, b_{01},b_{02},b_{03},b_{04}\right)
+ \left(c_{1i} | i \not \in \Lambda\right)
+ \left( b_{1i}-\alpha |i\in \Lambda \right), \cr
Q_{3\Lambda} &=
\C_0 + \left(s, b_{01},b_{02},b_{03},b_{04}\right)
+ \left(c_{1i} | i \not \in \Lambda\right)
+ \left( b_{1i}-b_{1j} |i,j\in \Lambda\right), \cr
Q_{4,2\alpha\beta} &=
 (s,c_{01},c_{04},c_{02}-c_{03},b_{01},b_{02},b_{03},b_{04},
 b_{11}-\alpha,b_{14}-\alpha,b_{12}-\beta,b_{13}-\beta) +\C_1. \cr
}
$$
These ideals are clearly prime ideals.
Let $x = f^3c_{21}b_{13}(b_{21}-b_{22})$ when $n > 2$
and $x = f^3$ otherwise.
We prove below that the $Q_{2\Lambda\alpha}$ and the $Q_{3\Lambda}$,
for all non-empty subsets $\Lambda \subseteq \{1,2,3,4\}$,
are associated primes of $J$,
and that when $n = 2$,
also the $Q_{4,2\alpha\beta}$ are associated.
Furthermore,
we prove that these $15d + d$ ($+ d^2$) primes
are the only new embedded primes of $J$ which do not contain $x$.

Consider the ideal
$$
\eqalignno{
\hat J &=(c_{01}-c_{02}b_{02}^d,c_{04}-c_{03}b_{03}^d)
+s(c_{02}-c_{03}) \cr
&\hskip2em+
c_{02}b_{02}^d (b_{02}^d-b_{01}^d,c_{02}-c_{03},b_{01}-b_{04},
b_{03}^d-b_{01}^d) \cr
&\hskip2em+
(c_{02}b_{01}-c_{03}b_{04},c_{02}(s-fb_{02}^d),c_{02}c_{1i}(b_{02}-b_{1i}b_{03}),
c_{02}b_{02}^d-c_{03}b_{03}^d) \cr
&\hskip2em+
c_{02}b_{02}^{2d}(c_{11}-c_{12},c_{14}-c_{13},c_{13}-c_{12},
c_{11}(b_{11}-b_{14}),c_{11}(b_{12}-b_{13}),
c_{11} \delta_{n \ge 3}).\cr
}
$$
It is easy to see that $x$ multiplies $\hat J$ into $J$
and that $\hat J$ contains $J$.
We will find all the associated primes of $\hat J$,
of which the only new ones are the $Q_{2\Lambda\alpha}$,
$Q_{3\Lambda}$, and $Q_{4,2\alpha\beta}$.
We will show that $x$ is not in any of the associated primes,
so that then $\hat J = J : x$.
Thus $\Ass(R/\hat J) = \Ass(R/J : x)$,
and every associated prime of $\hat J$ is also associated to $J$.
Thus it suffices to find all the associated primes of $\hat J$
and to show that $x$ is in none of them.

By Fact~\factassocses,
$\Ass \left( {R \over \hat J} \right)
\subseteq
 \Ass \left( {R \over \hat J : c_{02}} \right)
\cup \Ass \left( {R \over \hat J + (c_{02})} \right)$.
Note that
$\hat J + (c_{02}) =
(c_{01},c_{02},c_{04},$
$c_{03}b_{03}^d,sc_{03},c_{03}b_{04})$
$=p_0 \cap p_{-2}$ is an intersection of some minimal components of $J$,
and $x$ is a non-zerodivisor modulo each of them.
Hence it suffices to find the associated primes of $\hat J : c_{02}$,
and to show that $x$ is not in any of them:
$$
\eqalignno{
\hat J: c_{02} & =
(c_{01}-c_{02}b_{02}^d,c_{04}-c_{02}b_{02}^d,s-fb_{02}^d)\cr
&\hskip2em+
b_{02}^d (b_{02}^d-b_{01}^d,c_{02}-c_{03},b_{01}-b_{04},b_{03}^d-b_{01}^d)
+(c_{1i}(b_{02}-b_{1i}b_{03})) \cr
&\hskip2em+
b_{02}^{2d}(c_{11}-c_{12},c_{14}-c_{13},c_{13}-c_{12},
c_{11}(b_{11}-b_{14}),c_{11}(b_{12}-b_{13}),
c_{11} \delta_{n \ge 3})\cr
&\hskip2em+
(s(c_{02}-c_{03}),c_{02}b_{02}^d-c_{03}b_{03}^d,
c_{02}b_{01}-c_{03}b_{04}) : c_{02} \cr
&= (c_{01}-c_{02}b_{02}^d,c_{04}-c_{02}b_{02}^d,s-fb_{02}^d)
\cr
&\hskip2em+
b_{02}^d (b_{02}^d-b_{01}^d,c_{02}-c_{03},b_{01}-b_{04},b_{03}^d-b_{01}^d)
+(c_{1i}(b_{02}-b_{1i}b_{03})) \cr
&\hskip2em+
b_{02}^{2d}(c_{11}-c_{12},c_{14}-c_{13},c_{13}-c_{12},
c_{11}(b_{11}-b_{14}),c_{11}(b_{12}-b_{13}),
c_{11} \delta_{n \ge 3})\cr
&\hskip2em+
(c_{02}b_{02}^d-c_{03}b_{03}^d,
c_{02}b_{01}-c_{03}b_{04},
b_{01}b_{03}^d-b_{04}b_{02}^d) \cr
&\hskip2em
+ s(c_{02}-c_{03},b_{01}-b_{04},b_{02}^d-b_{03}^d) \cr
&= (c_{01}-c_{02}b_{02}^d,c_{04}-c_{02}b_{02}^d,s-fb_{02}^d)
\cr
&\hskip2em+
b_{02}^d (b_{02}^d-b_{01}^d,c_{02}-c_{03},b_{01}-b_{04},b_{03}^d-b_{01}^d)
 +(c_{1i}(b_{02}-b_{1i}b_{03})) \cr
&\hskip2em+
b_{02}^{2d}(c_{11}-c_{12},c_{14}-c_{13},c_{13}-c_{12},
c_{11}(b_{11}-b_{14}),c_{11}(b_{12}-b_{13}),
c_{11} \delta_{n \ge 3})\cr
&\hskip2em+
(c_{02}b_{02}^d-c_{03}b_{03}^d,
c_{02}b_{01}-c_{03}b_{04},
b_{01}b_{03}^d-b_{04}b_{02}^d). \cr
}
$$
Again by Fact~\factassocses,
$\Ass \left( {R \over \hat J : c_{02}} \right)
\subseteq
 \Ass \left( {R \over \hat J : c_{02}b_{02}^d} \right)
\cup \Ass \left( {R \over (\hat J:c_{02}) + (b_{02}^d)} \right)$.
Note that
$$
(\hat J: c_{02}) + (b_{02}^d) =
(c_{01},c_{04},s,b_{02}^d,c_{1i}(b_{02}-b_{1i}b_{03}),c_{03}b_{03}^d,
c_{02}b_{01}-c_{03}b_{04}, b_{01}b_{03}^d),
$$
which decomposes as
$$
\eqalignno{
&=
(c_{01},c_{04},s,b_{02}^d,c_{1i}(b_{02}-b_{1i}b_{03}),c_{1i}b_{1i}^d,
c_{1i}c_{1j}(b_{1i}-b_{1j}),
c_{03}, b_{01}) \cr
&\hskip1em\cap
(c_{01},c_{04},s,b_{02}^d,c_{1i}(b_{02}-b_{1i}b_{03}),b_{03}^d,
c_{02}b_{01}-c_{03}b_{04}) \cr
&= p_{-4} \cap p_{-3} \cap q_1, \cr
}
$$
as in  Section~\sectsixteen.
Thus $x$ is a non-zerodivisor on this ideal,
and no new associated primes appear.
Thus it suffices to find the associated primes of $\hat J :c_{02}b_{02}^d$:
$$
\eqalignno{
\hat J &: c_{02}b_{02}^d
= (c_{01}-c_{02}b_{02}^d,c_{04}-c_{02}b_{02}^d,s-fb_{02}^d)
\cr
&\hskip2em+
(b_{02}^d-b_{01}^d,c_{02}-c_{03},b_{01}-b_{04},
b_{03}^d-b_{01}^d) \cr
&\hskip2em+
b_{02}^d (c_{11}-c_{12},c_{14}-c_{13},c_{13}-c_{12},
c_{11}(b_{11}-b_{14}),c_{11}(b_{12}-b_{13}),
c_{11} \delta_{n \ge 3})\cr
&\hskip2em+
(c_{1i}(b_{02}-b_{1i}b_{03}),
c_{02}b_{02}^d-c_{03}b_{03}^d,
c_{02}b_{01}-c_{03}b_{04},
b_{01}b_{03}^d-b_{04}b_{02}^d) : b_{02}^d. \cr
}
$$
The next two displays will compute the colon ideal in the last row.
As in the computation of $p_{-4}$ in [S2],
\pagelabel{\qcalcc}
$$
\eqalignno{
&(c_{1i}(b_{02}-b_{1i}b_{03}),
c_{02}b_{02}^d-c_{03}b_{03}^d,
c_{02}b_{01}-c_{03}b_{04},
b_{01}b_{03}^d-b_{04}b_{02}^d) \cr
&= \bigcap_{\Lambda}
\Bigl((c_{02}b_{02}^d-c_{03}b_{03}^d,
c_{02}b_{01}-c_{03}b_{04},
b_{01}b_{03}^d-b_{04}b_{02}^d,
c_{1i},b_{02}-b_{1j}b_{03} | i \not \in \Lambda, j \in \Lambda)
\Bigr) \cr
&=\Bigl((c_{02}b_{02}^d-c_{03}b_{03}^d,
c_{02}b_{01}-c_{03}b_{04}, b_{01}b_{03}^d-b_{04}b_{02}^d)
+ \C_1 \Bigr) \cr
&\bigcap_{\Lambda \not = \emptyset}
\Bigl(((c_{02}b_{1j}^d-c_{03})b_{03}^d,
c_{02}b_{01}-c_{03}b_{04},
(b_{01}-b_{04}b_{1j}^d)b_{03}^d,
c_{1i},b_{02}-b_{1j}b_{03} | i \not \in \Lambda, j \in \Lambda)
\Bigr), \cr
}
$$
which coloned with $b_{02}^d$ equals
$$
\eqalignno{
&\Bigl((c_{02}b_{02}^d-c_{03}b_{03}^d,
c_{02}b_{01}-c_{03}b_{04}, b_{01}b_{03}^d-b_{04}b_{02}^d)
+ \C_1 \Bigr) \cr
&\bigcap_{\Lambda \not = \emptyset}
\Bigl((c_{02}b_{1j}^d-c_{03},
b_{01}-b_{04}b_{1j}^d,
c_{1i},b_{02}-b_{1j}b_{03}, b_{1j}-b_{1j'} | i \not \in \Lambda, j,j' \in \Lambda)
\Bigr) \cr
&=\left(
c_{02}b_{02}^d-c_{03}b_{03}^d,
c_{02}b_{01}-c_{03}b_{04},
b_{01}b_{03}^d-b_{04}b_{02}^d,
c_{1i}(b_{02}-b_{1i}b_{03})\right) \cr
&\hskip2em +\left(
c_{1i}c_{1j}(b_{1i}-b_{1j}),
c_{1i}(b_{01}-b_{1i}^db_{04}),
c_{1i}(c_{03}-b_{1i}^dc_{02}) | i,j \right). \cr
}
$$
Thus
$$
\eqalignno{
\hat J &: c_{02}b_{02}^d
= (c_{01}-c_{02}b_{02}^d,c_{04}-c_{02}b_{02}^d,s-fb_{02}^d, c_{02}-c_{03},
b_{02}^d-b_{01}^d,b_{01}-b_{04},b_{03}^d-b_{01}^d) \cr
&\hskip2em+
b_{02}^d (c_{11}-c_{12},c_{14}-c_{13},c_{13}-c_{12},
c_{11}(b_{11}-b_{14}),c_{11}(b_{12}-b_{13}),
c_{11} \delta_{n \ge 3})\cr
&\hskip2em +\left(
c_{1i}(b_{02}-b_{1i}b_{03}),
c_{1i}c_{1j}(b_{1i}-b_{1j}),
c_{1i}b_{01}(1-b_{1i}^d),
c_{02}c_{1i}(1-b_{1i}^d) | i,j \right). \cr
}
$$
Now let $J' = \hat J : c_{02}b_{02}^{2d}$
and $J'' = (\hat J : c_{02}b_{02}^d) + (b_{02}^d)$.
By Fact~\factassocses,
the set of associated primes of $\hat J : c_{02}b_{02}^d$
is contained in the union of the sets of associated primes of $J'$
and $J''$.

First we analyze $J''$:
$$
\eqalignno{
J'' &= (c_{01},c_{04},c_{02}-c_{03},s, b_{0i}^d,b_{01}-b_{04}) \cr
&\hskip2em +\left(
c_{1i}(b_{02}-b_{1i}b_{03}),
c_{1i}c_{1j}(b_{1i}-b_{1j}),
c_{1i}b_{01}(1-b_{1i}^d),
c_{02}c_{1i}(1-b_{1i}^d) | i,j \right). \cr
}
$$
This decomposes as follows:
$$
\eqalignno{
J''&= (\C_0 +\left( s, b_{0i}^d,b_{01}-b_{04},c_{1i}(b_{02}-b_{1i}b_{03}),
c_{1i}c_{1j}(b_{1i}-b_{1j}),
c_{1i}b_{01}(1-b_{1i}^d)\right)) \cr
&\hskip1em\cap \left(c_{01},c_{04},c_{02}-c_{03},s, b_{0i}^d,b_{01}-b_{04},
c_{1i}(b_{02}-b_{1i}b_{03}), c_{1i}c_{1j}(b_{1i}-b_{1j}),
c_{1i}(1-b_{1i}^d) \right). \cr
}
$$
Clearly the ideal $q_2$ in the second row decomposes
\pagelabel{\pageminussix}
as the intersection of $Q_{2\Lambda\alpha}$-primary components,
and the ideal in the first row decomposes as
$$
\eqalignno{
&(\C_0 + \left(s, b_{0i}^d,b_{01}-b_{04},
c_{1i}(b_{02}-b_{1i}b_{03}),
c_{1i}c_{1j}(b_{1i}-b_{1j}),
c_{1i}(1-b_{1i}^d)\right)) \cr
&\hskip1em\cap (\C_0 + \left(s, b_{0i}^d,b_{01},b_{04},
c_{1i}(b_{02}-b_{1i}b_{03}),
c_{1i}c_{1j}(b_{1i}-b_{1j})\right)). \cr
}
$$
Of these components,
the ideal $q_3$ in the second row
is an intersection of the $Q_{3\Lambda}$-primary components,
and the ideal in the first row contains $q_2$,
and is thus redundant for computing the associated primes of $J''$.
Thus the set of associated primes of $J''$ is a subset of
$\{Q_{2\Lambda\alpha}, Q_{3\Lambda}\}$.
Clearly $x$ is not in any $Q_{2\Lambda\alpha}$ and $Q_{3\Lambda}$.

It remains to compute a decomposition of $J'$:
$$
\eqalignno{
J' &= (c_{01}-c_{02}b_{02}^d,c_{04}-c_{02}b_{02}^d,s-fb_{02}^d,c_{02}-c_{03},
b_{02}^d-b_{01}^d,b_{01}-b_{04},b_{03}^d-b_{01}^d) \cr
&\hskip2em+
(c_{11}-c_{12},c_{14}-c_{13},c_{13}-c_{12},
c_{11}(b_{11}-b_{14}),c_{11}(b_{12}-b_{13}),
c_{11} \delta_{n \ge 3})\cr
&\hskip2em +\left(
c_{11}(b_{02}-b_{1i}b_{03}),
c_{11}^2(b_{1i}-b_{1j}),
c_{11}(1-b_{1i}^d)| i,j \right). \cr
}
$$
By coloning and adding $c_{11}^2$:
$$
\eqalignno{
J' &=(c_{01}-c_{02}b_{02}^d,c_{04}-c_{02}b_{02}^d,s-fb_{02}^d,c_{02}-c_{03},
b_{01}-b_{04},b_{03}^d-b_{01}^d) \cr
&\hskip2em+
\left(c_{11}-c_{12},c_{14}-c_{13},c_{13}-c_{12},\delta_{n \ge 3},
b_{02}-b_{1i}b_{03},
b_{1i}-b_{1j}, 1-b_{1i}^d\right)) \cr
&\hskip1em\cap ((c_{01}-c_{02}b_{02}^d,c_{04}-c_{02}b_{02}^d,s-fb_{02}^d,
c_{02}-c_{03}, b_{02}^d-b_{01}^d,b_{01}-b_{04},b_{03}^d-b_{01}^d) \cr
&\hskip2em+
(c_{11}-c_{12},c_{14}-c_{13},c_{13}-c_{12},
c_{11}(b_{11}-b_{14}),c_{11}(b_{12}-b_{13}),
c_{11} \delta_{n \ge 3})\cr
&\hskip2em +\left(
c_{11}(b_{02}-b_{1i}b_{03}),
c_{11}^2,
c_{11}(1-b_{1i}^d)\right)) \cr
&= p_2 \delta_{n=2} \cr
&\hskip1em\cap
((c_{01}-c_{02}b_{02}^d,c_{04}-c_{02}b_{02}^d,s-fb_{02}^d,c_{02}-c_{03},
b_{02}^d-b_{01}^d,b_{01}-b_{04},b_{03}^d-b_{01}^d) + \C_1) \cr
&\hskip1em\cap ((c_{01}-c_{02}b_{02}^d,c_{04}-c_{02}b_{02}^d,s-fb_{02}^d,
c_{02}-c_{03}, b_{02}^d-b_{01}^d,b_{01}-b_{04},b_{03}^d-b_{01}^d) + \D_1 \cr
&\hskip2em+
\left(b_{11}-b_{14},b_{12}-b_{13},\delta_{n \ge 3},b_{02}-b_{1i}b_{03},
c_{11}^2,1-b_{1i}^d\right)). \cr
}
$$
By coloning and adding $b_{03}$ on the third component,
$$
\eqalignno{
J' &=
p_2 \delta_{n=2} \cap p_1 \cr
&\cap ((c_{01}-c_{02}b_{02}^d,c_{04}-c_{02}b_{02}^d,s-fb_{02}^d,
c_{02}-c_{03}, b_{02}^d-b_{01}^d,b_{01}-b_{04},b_{03}^d-b_{01}^d) + \D_1 \cr
&\hskip1em+
\left(\delta_{n \ge 3},b_{02}-b_{1i}b_{03},
b_{1i}-b_{1j}, c_{11}^2,1-b_{1i}^d\right)). \cr
&\cap ((c_{01},c_{04},s,c_{02}-c_{03},
b_{01}^d,b_{01}-b_{04},b_{02},b_{03},
b_{11}-b_{14},b_{12}-b_{13},\delta_{n \ge 3}, c_{11}^2,1-b_{1i}^d) + \D_1). \cr
}
$$
The second to the last ideal above properly contains $p_2 \delta_{n=2}$,
and the last ideal $q_{4,2}$
is an intersection of $Q_{4,2\alpha\beta}$-primary components
when $n = 2$.
This proves that
$$
J' = \hat J : c_{02}b_{02}^{2d} =
p_1 \cap p_2 \delta_{n=2} \cap q_{4,2}.
$$
As $x$ is a non-zerodivisor modulo this ideal,
this also finishes the proof that $\hat J = J : x$.
Furthermore,
this proves that the set of new embedded primes of
$J$ which do not contain $x$
is contained in the set of associated primes of the ideal
$\hat J : c_{02}b_{02}^d$,
and that this latter set is a subset of
$$
\{Q_{2\Lambda\alpha}, Q_{3\Lambda}, Q_{4,2\alpha\beta}\}.
$$
It remains to prove that the prime ideals
$Q_{2\emptyset\alpha}$ are not associated to $J$,
and that every element of
$$
\{Q_{2\Lambda\alpha}, Q_{3\Lambda}, Q_{4,2\alpha\beta}\delta_{n=2}
| \Lambda \not = \emptyset\}
$$
is associated to $J$.
Clearly,
when $n = 2$,
as $J' : c_{02} = p_2 \cap q_{4,2} \not = p_2$,
all the $Q_{4,2\alpha\beta}$ are associated to $J'$ and thus to $J$.

Let $\Lambda$ be a subset of $\{1,2,3,4\}$.
Let $K$ be the ideal $\hat J : c_{02}b_{02}^d$ coloned with
a power of the element
$$
\prod_{i \in \Lambda \atop j \not \in \Lambda}
c_{1i}(1-b_{1j}^d)c_{02}.
$$
Note that this element is contained in $Q_{4,2\alpha\beta}$,
in $Q_{3\Lambda}$,
in $p_1$ if $\Lambda \not = \emptyset$,
and in $p_2$ if $\Lambda \not = \{1,2,3,4\}$.
If $\Lambda = \emptyset$,
then $K = p_1$,
so $Q_{2\Lambda\alpha}$ are all redundant.
If however, $\Lambda \not = \emptyset$,
and if $\Lambda = \{1,2,3,4\}$,
then $K \not = p_2\delta_{n=2}$,
which proves that the $Q_{2\Lambda\alpha}$ are all associated to $J$.

Similarly,
by coloning with
$\prod_{i \in \Lambda \atop j,j' \not \in \Lambda}
c_{1i}(b_{1j}-b_{1j'})(b_{1i}-b_{1j})\prod_{j=1}^4(1-b_{1j}^d)$
we get that each $Q_{3\Lambda}$ is associated to $J$
if and only if $\Lambda \not = \emptyset$.
Thus

\thm
\label{\thmnox}
Set $x =fc_{21}b_{13}(b_{21}-b_{22})$ when $n > 2$
and $x = f$ otherwise.
Then the set of embedded primes of $J$ which do not contain $x$
is contained in
$$
\lbrace Q_{1\Lambda},
Q_{2\Lambda'\alpha},
Q_{3\Lambda'},
Q_{4,2\alpha\beta}\delta_{n=2}| \Lambda' \not = \emptyset\rbrace,
$$
and each of these primes is associated to $J$.
\qed
\endb

For clarity we record these embedded primes in a table:

$$
\displayindent =1em
\halign{\mvrule \hskip 0.1em #\hfil \hskip 0.1em
&& \mvrule \hskip 0.5em \relax #\hfil \mvrule \cr
\noalign{\hrule}
embedded prime ($\Lambda \not = \emptyset, \alpha^d = 1, \beta^d = 1$)
& height  \cr
\noalign{\hrule}
\noalign{\hrule}
$Q_{2\Lambda\alpha} =
\left(s, c_{01}, c_{03}-c_{02}, c_{04}, b_{01},b_{02},b_{03},b_{04}\right)
+ \left(c_{1i} | i \not \in \Lambda\right)
+ \left( b_{1i}-\alpha |i\in \Lambda \right)$
& 12 \cr
$Q_{3\Lambda} =
\C_0 + \left(s, b_{01},b_{02},b_{03},b_{04}\right)
+ \left(c_{1i} | i \not \in \Lambda\right)
+ \left( b_{1i}-b_{1j} |i,j\in \Lambda\right)$
& 12 \cr
$Q_{4,2\alpha\beta} =
 (s,c_{01},c_{04},c_{02}-c_{03},b_{01},b_{02},b_{03},b_{04})$ & 16 \cr
\hskip 4em$(b_{11}-\alpha,b_{14}-\alpha,b_{12}-\beta,b_{13}-\beta)+\C_1$
\hfill if $n = 2$ only & \cr
\noalign{\hrule}
}
$$
\vskip2ex

%
%

\vfill\eject
\section{$\bf (n-1)(d^3-d)$ more embedded primes, for $\bf n > 2$}
\sectlabel{\sectmoremoreemb}
\pagelabel{\qsectmoremoreemb}

The embedded primes of $J$ found so far do not contain $b_{2i}-b_{2j}$.
Without this assumption there are many more embedded primes of $J$,
and their number grows with $n$ and $d$.
In this section,
$(n-1)(d^3-d)$ more embedded primes are found in the case when $n > 2$.
The main theorem of this section,
Theorem~\thmthirdemb,
says that these primes are the only new ones not containing the element $x$,
where $x$ is
$$
x = \left\lbrace
\matrix{f^3(c_{21}\cdots c_{r-1,1})
b_{13}^{2d+1}(b_{23}\cdots b_{r-1,3})(1-b_{r1}), \hfill&
\hbox{\ \ if $r <n$}, \hfill\cr
f^3(c_{21}\cdots c_{r-1,1})b_{13}^{2d+1}(b_{23}\cdots b_{r-1,3}), \hfill&
\hbox{\ \ otherwise}.\hfill \cr
}
\right.
$$
Throughout this section,
$n > 2$.

For each $r \in \{2,\ldots, n\}$
and $\alpha, \beta$ and $\gamma$ in $k$
such that $\alpha^d = \beta^d = \gamma^d = 1$,
define
$$
\eqalignno{
&Q_{4r\alpha\beta\gamma} =
 (s,c_{01},c_{04},c_{02}-c_{03},b_{01},b_{02},b_{03},b_{04}) \cr
&\hskip2em+ (b_{12}-b_{2i}b_{13},b_{2i}-b_{2j}) \delta_{r > 2}
 + (b_{11}-\alpha,b_{14}-\alpha,b_{12}-\beta,b_{13}-\gamma) \cr
&\hskip2em+\C_1 +\D_2 + \cdots + \D_{r-1} + \C_r + B_{3,r-1}. \cr
}
$$
It is proved in this section that these prime ideals are associated to $J$
if and only if $\{\alpha, \beta, \gamma\}| > 1$,
and that these $(n-1)(d^3-d)$ prime ideals
are the only new associated prime ideals of $J$
which do not contain the element $x$ defined above.
 
For all $2 \le r \le n$,
with the convention that $c_{ni} = b_{ni} = 1$, $\C_n = (0)$,
all these cases can be analyzed simultaneously.
Consider the ideal
$$
\eqalignno{
K &=\left(c_{01}-c_{02}b_{02}^d,c_{01}-c_{04},
s(c_{02}-c_{03}),c_{02}b_{01}-c_{03}b_{04},c_{02}(s-fb_{02}^d),
c_{02}b_{02}^d-c_{03}b_{03}^d\right) \cr
&\hskip2em+
c_{02}(b_{02}^d,c_{13}b_{03}^d) (b_{02}^d-b_{01}^d,c_{02}-c_{03},b_{01}-b_{04},
b_{03}^d-b_{01}^d) +
(c_{02}c_{1i}(b_{02}-b_{1i}b_{03}))  \cr
&\hskip2em+c_{02}b_{02}^{2d}(c_{1i}b_{1i}^d - c_{13}b_{13}^d)
+c_{02}c_{13}b_{03}^{2d}\left( \D_1 + (b_{11}-b_{14},1-b_{1i}^d)\right) \cr
&\hskip2em
+c_{02}c_{13}b_{03}^{2d}\left(
(c_{11},b_{02},b_{03})(1-b_{2i})
+\left(b_{12}-b_{2i}b_{13},b_{2i}-b_{2j}\right)\right)\delta_{r>2} \cr
&\hskip2em+
\sum_{k=2}^{r-2} c_{02}b_{03}^{2d}c_{13}
\left( \D_k +(1-b_{k+1,i})\right) +
 c_{02}c_{13}b_{03}^{2d}\left(\D_{r-1} + \C_r\right). \cr
}
$$
It is easy to see that $K$ contains $J$ and
that $x$ multiplies $K$ into $J$,
except possibly that $x$ multiplies
the element $c_{02}b_{03}^{2d}c_{11}c_{13}(1-b_{2i})$ into $J$:
$$
\eqalignno{
c_{02}b_{03}^{2d}c_{11}c_{13}(1-b_{2i})x &\in
(f^2c_{02}b_{03}^{2d}c_{11}c_{13}(1-b_{2i})c_{21}b_{13}^{2d+1}) + J \cr
&=(f^2c_{02}b_{02}^{2d}c_{11}c_{13}(1-b_{2i})c_{21}b_{13}) + J \cr
&=(sfc_{01}c_{11}^2(1-b_{2i})c_{2i}b_{13}) + J \cr
&=(sfc_{04}c_{11}c_{12}(1-b_{2i})c_{2i}b_{13}) + J \cr
&=(sfc_{04}c_{11}c_{12}(b_{13}-b_{12})c_{2i}) + J \cr
&=(sfc_{01}c_{11}(c_{13}b_{13}-c_{12}b_{12})c_{2i}) + J \cr
&=(sfc_{02}c_{11}(c_{13}b_{13}-c_{12}b_{12})b_{03}^d c_{2i}) + J \cr
&=(sfc_{02}c_{11}(c_{13}-c_{12})b_{02}b_{03}^{d-1} c_{2i}) + J \cr
&=(sfc_{02}c_{11}b_{11}(c_{13}-c_{12})b_{03}^d c_{2i}) + J \cr
&=(sfc_{01}c_{11}b_{11}(c_{13}-c_{12}) c_{2i}) + J = J. \cr
}
$$
The intermediate goal in this section is to find a primary decomposition of $K$.
It turns out that $x$ is a non-zerodivisor on $K$,
which proves that $K = J : x$,
and thus determines all the associated primes of $J$ which do not contain $x$.

\vskip 3ex
By Fact~\factassocses,
$\Ass \left( {R \over K} \right)
\subseteq
 \Ass \left( {R \over K :c_{02}} \right)
\cup \Ass \left( {R \over K + (c_{02})} \right)$.
The second set is easy:
$K + (c_{02})$ equals
$(c_{01},c_{02},c_{04}, sc_{03},c_{03}b_{04}, c_{03}b_{03}^d)$
$=p_0 \cap p_{-2}$,
which is an intersection of some minimal components of $J$
(none of which contain $x$),
so it suffices to find the associated primes of $K : c_{02}$.
By Facts~\factcolony\ and \factfewvars:
$$
\eqalignno{
&K : c_{02}  =
\left(c_{01}-c_{02}b_{02}^d,c_{01}-c_{04},s-fb_{02}^d\right) \cr
&\hskip2em+
(b_{02}^d,c_{13}b_{03}^d) (b_{02}^d-b_{01}^d,c_{02}-c_{03},b_{01}-b_{04},
b_{03}^d-b_{01}^d) +
(c_{1i}(b_{02}-b_{1i}b_{03}))  \cr
&\hskip2em+b_{02}^{2d}(c_{1i}b_{1i}^d - c_{13}b_{13}^d)
+c_{13}b_{03}^{2d}\left( \D_1 + (b_{11}-b_{14},1-b_{1i}^d)\right) \cr
&\hskip2em
+c_{13}b_{03}^{2d}\left(
(c_{11},b_{02},b_{03})(1-b_{2i})
+\left(b_{12}-b_{2i}b_{13},b_{2i}-b_{2j}\right)\right)\delta_{r>2} \cr
&\hskip2em+
\sum_{k=2}^{r-2} b_{03}^{2d}c_{13} \left( \D_k +(1-b_{k+1,i})\right)
+ c_{13}b_{03}^{2d}\left(\D_{r-1} + \C_r\right) \cr
&\hskip2em+
(s(c_{02}-c_{03}),
c_{02}b_{01}-c_{03}b_{04},c_{02}b_{02}^d-c_{03}b_{03}^d) : c_{02}. \cr
}
$$
The latter colon ideal equals
$$
(c_{02}b_{02}^d-c_{03}b_{03}^d,
c_{02}b_{01}-c_{03}b_{04},
b_{01}b_{03}^d-b_{04}b_{02}^d)
+ s(c_{02}-c_{03},b_{01}-b_{04},b_{02}^d-b_{03}^d),
$$
so that
$$
\eqalignno{
&K : c_{02}  =
\left(c_{01}-c_{02}b_{02}^d,c_{01}-c_{04},s-fb_{02}^d\right) \cr
&\hskip2em+
(b_{02}^d,c_{13}b_{03}^d) (b_{02}^d-b_{01}^d,c_{02}-c_{03},b_{01}-b_{04},
b_{03}^d-b_{01}^d) +
(c_{1i}(b_{02}-b_{1i}b_{03}))  \cr
&\hskip2em+b_{02}^{2d}(c_{1i}b_{1i}^d - c_{13}b_{13}^d)
+c_{13}b_{03}^{2d}\left( \D_1 + (b_{11}-b_{14},1-b_{1i}^d)\right) \cr
&\hskip2em
+c_{13}b_{03}^{2d}\left(
(c_{11},b_{02},b_{03})(1-b_{2i})
+\left(b_{12}-b_{2i}b_{13},b_{2i}-b_{2j}\right)\right)\delta_{r>2} \cr
&\hskip2em+
\sum_{k=2}^{r-2} b_{03}^{2d}c_{13} \left( \D_k +(1-b_{k+1,i})\right)
+ c_{13}b_{03}^{2d}\left(\D_{r-1} + \C_r\right) \cr
&\hskip2em+
(c_{02}b_{02}^d-c_{03}b_{03}^d,
c_{02}b_{01}-c_{03}b_{04},
b_{01}b_{03}^d-b_{04}b_{02}^d). \cr
}
$$
\pagelabel{\qstep}
Again by Fact~\factassocses,
$\Ass \left( {R \over K : c_{02}} \right)
\subseteq
 \Ass \left( {R \over K :c_{02}b_{03}^d} \right)
\cup \Ass \left( {R \over (K :c_{02}) + (b_{03}^d)} \right)$.
Note that
$$
\eqalignno{
&(K :c_{02}) + (b_{03}^d)  =
(c_{01}-c_{02}b_{02}^d,c_{01}-c_{04},
s-fb_{02}^d,b_{03}^d)+(c_{1i}(b_{02}-b_{1i}b_{03})) \cr
&\hskip2em+
b_{02}^d(b_{01}^d,b_{02}^d,c_{02}-c_{03},b_{01}-b_{04})
+(c_{02}b_{02}^d,c_{02}b_{01}-c_{03}b_{04},b_{04}b_{02}^d). \cr
}
$$
By Fact~\factassocses,
$\Ass \left( {(K :c_{02}) + (b_{03}^d)} \right)
\subseteq
 \Ass \left( {R \over ((K :c_{02}) + (b_{03}^d)):b_{02}^d} \right)
\cup \Ass \left( {R \over (K :c_{02}) + (b_{03}^d,b_{02}^d)} \right)$.
Then
$$
\eqalignno{
&(K :c_{02}) + (b_{03}^d,b_{02}^d)  =
(s,c_{01},c_{04},b_{02}^d,b_{03}^d,c_{1i}(b_{02}-b_{1i}b_{03}),
c_{02}b_{01}-c_{03}b_{04}) \cr
&=
(s,c_{01},c_{04},b_{02},b_{03},c_{02}b_{01}-c_{03}b_{04}) \cr
&\hskip1em\cap
(s,c_{01},c_{04},b_{02}^d,b_{03}^d,c_{1i}(b_{02}-b_{1i}b_{03}),
c_{1i}c_{1j}(b_{1i}-b_{1j}),c_{02}b_{01}-c_{03}b_{04}). \cr
}
$$
The first component is $P_{-3}$,
and the second component is the intersection of ideals
primary to the $Q_{1\Lambda}$,
as $\Lambda$ varies over the subsets of $\{1,2,3,4\}$.
None of these prime ideals contains $x$.

Next,
$((K :c_{02}) + (b_{03}^d)):b_{02}^d$ equals
$$
\eqalignno{
&=\C_0 + (s,b_{01},b_{02}^d,b_{03}^d,b_{04})
+(b_{03}^d,c_{1i}(b_{02}-b_{1i}b_{03}), c_{02}b_{01}-c_{03}b_{04}) :b_{02}^d\cr
&=
\C_0 + (s,b_{01},b_{02}^d,b_{03}^d,b_{04})+\C_1, \cr
}
$$
and again $x$ is a non-zerodivisor modulo this ideal,
and the associated prime of this ideal ($Q_{3\emptyset}$)
is not associated to $J$ by Theorem~\thmnox.

This finishes the analysis of the associated primes of $(K : c_{02}) + (b_{03}^d)$.
It remains to analyze $K : c_{02}b_{03}^d$.
This colon ideal is
$$
\eqalignno{
&\hskip1em K : c_{02} b_{03}^d =
\left(c_{01}-c_{02}b_{02}^d,c_{01}-c_{04},s-fb_{02}^d\right) \cr
&\hskip2em+
c_{13}(b_{02}^d-b_{01}^d,c_{02}-c_{03},b_{01}-b_{04}, b_{03}^d-b_{01}^d)
\cr
&\hskip2em+b_{02}^d (c_{1i}b_{1i}^d - c_{13}b_{13}^d)
+c_{13}b_{03}^d \left( \D_1 + (b_{11}-b_{14},1-b_{1i}^d)\right) \cr
&\hskip2em
+c_{13}b_{03}^d \left(
(c_{11},b_{02},b_{03})(1-b_{2i})
+\left(b_{12}-b_{2i}b_{13},b_{2i}-b_{2j}\right)\right)\delta_{r>2} \cr
&\hskip2em+
\sum_{k=2}^{r-2} b_{03}^d c_{13} \left( \D_k +(1-b_{k+1,i})\right)
+ c_{13}b_{03}^d \left(\D_{r-1} + \C_r\right) +
\left(L : b_{03}^d\right), \cr
}
$$
where
$$
\eqalignno{
L &=
(c_{1i}(b_{02}-b_{1i}b_{03})) +
b_{02}^d(b_{02}^d-b_{01}^d,c_{02}-c_{03},b_{01}-b_{04},b_{03}^d-b_{01}^d)
 \cr
&\hskip2em+
(c_{02}b_{02}^d-c_{03}b_{03}^d,
c_{02}b_{01}-c_{03}b_{04},
b_{01}b_{03}^d-b_{04}b_{02}^d). \cr
}
$$
The next two pages will compute $L : b_{03}^d$.
First of all,
coloning with $b_{02}^d$ gives:
$$
\eqalignno{
L : b_{02}^d &=
(b_{02}^d-b_{01}^d,c_{02}-c_{03},b_{01}-b_{04},b_{03}^d-b_{01}^d) \cr
&\hskip2em+\left(c_{1i}(b_{02}-b_{1i}b_{03}),
c_{02}b_{02}^d-c_{03}b_{03}^d,
c_{02}b_{01}-c_{03}b_{04},
b_{01}b_{03}^d-b_{04}b_{02}^d\right) : b_{02}^d, \cr
}
$$
which by a computation on page~\qcalcc\ equals
$$
\eqalignno{
L : b_{02}^d &=
(b_{02}^d-b_{01}^d,c_{02}-c_{03},b_{01}-b_{04},b_{03}^d-b_{01}^d) \cr
&\hskip2em+
(c_{1i}(b_{02}-b_{1i}b_{03}),
c_{1i}c_{1j}(b_{1i}-b_{1j}),
b_{01}c_{1i}(1-b_{1i}^d),
c_{02}c_{1i}(1-b_{1i}^d)), \cr
}
$$
so that $L : b_{02}^d b_{01}$ equals
$$
(b_{02}^d-b_{01}^d,c_{02}-c_{03},b_{01}-b_{04},b_{03}^d-b_{01}^d,
c_{1i}(b_{02}-b_{1i}b_{03}),
c_{1i}c_{1j}(b_{1i}-b_{1j}),
c_{1i}(1-b_{1i}^d)).
$$
Note that neither $b_{01}$ nor $b_{02}$ is a zero-divisor modulo
$L : b_{02}^d b_{01}$,
so that by Fact~\factdec,
$$
\eqalignno{
&L= (L : b_{02}^d b_{01}) \cap (L + (b_{02}^d b_{01})) \cr
&=
\Bigl(b_{02}^d-b_{01}^d,c_{02}-c_{03},b_{01}-b_{04},b_{03}^d-b_{01}^d,
c_{1i}(b_{02}-b_{1i}b_{03}),
c_{1i}c_{1j}(b_{1i}-b_{1j}),
c_{1i}(1-b_{1i}^d)\Bigr)\cr
&\hskip1em\cap
\Bigl((c_{1i}(b_{02}-b_{1i}b_{03})) +
b_{02}^d(b_{01},b_{02}^d,b_{03}^d,c_{02}-c_{03},b_{04})
 \cr
&\hskip2em+
(c_{02}b_{02}^d-c_{03}b_{03}^d,
c_{02}b_{01}-c_{03}b_{04},
b_{01}b_{03}^d)\Bigr) \cr
&=
\Bigl(b_{02}^d-b_{01}^d,c_{02}-c_{03},b_{01}-b_{04},b_{03}^d-b_{01}^d,
c_{1i}(b_{02}-b_{1i}b_{03}),
c_{1i}c_{1j}(b_{1i}-b_{1j}),
c_{1i}(1-b_{1i}^d)\Bigr)\cr
&\hskip1em\cap
\Bigl(\C_1 +
b_{02}^d(b_{01},b_{02}^d,b_{03}^d,c_{02}-c_{03},b_{04})
+(c_{02}b_{02}^d-c_{03}b_{03}^d,
c_{02}b_{01}-c_{03}b_{04},b_{01}b_{03}^d)\Bigr)\cr
&\hskip1em\bigcap_{\Lambda\not=\emptyset}
\Bigl((c_{02}b_{01}-c_{03}b_{04},b_{01}b_{03}^d) +(c_{1i}| i \not \in \Lambda) \cr
&\hskip2em+ (b_{02}-b_{1i}b_{03},
b_{1i}^db_{03}^d(b_{01},b_{03}^d,c_{02}-c_{03},b_{04}),
b_{03}^d(c_{02}b_{1i}^d-c_{03})
| i \in \Lambda)\Bigr). \cr
}
$$
This is still part of the effort to compute $L : b_{03}^d$.
Coloning the second component above with $b_{03}^d$ equals
$$
\eqalignno{
&\C_1 + (b_{01},b_{02}^d) +
\Bigl(
b_{02}^d(b_{01},b_{02}^d,c_{02}-c_{03},b_{04})
+(c_{02}b_{02}^d-c_{03}b_{03}^d,
c_{02}b_{01}-c_{03}b_{04})\Bigr) : b_{03}^d. \cr
}
$$
But
$$
\eqalignno{
&b_{02}^d(b_{01},b_{02}^d,c_{02}-c_{03},b_{04})
+(c_{02}b_{02}^d-c_{03}b_{03}^d,
c_{02}b_{01}-c_{03}b_{04}) \cr
&= \Bigl(b_{02}^{2d}, b_{01},c_{02}-c_{03},b_{04},
c_{02}b_{02}^d-c_{03}b_{03}^d\Bigr) \cap\Bigl(b_{02}^d,c_{03}b_{03}^d,
c_{02}b_{01}-c_{03}b_{04}\Bigr), \cr
}
$$
so that
$$
\eqalignno{
&\C_1 + (b_{01},b_{02}^d) +
\Bigl(
b_{02}^d(b_{01},b_{02}^d,c_{02}-c_{03},b_{04})
+(c_{02}b_{02}^d-c_{03}b_{03}^d,
c_{02}b_{01}-c_{03}b_{04})\Bigr) : b_{03}^d=\cr
&=\C_1 + (b_{01},b_{02}^d)
+ \Bigl(b_{02}^{2d}, b_{01},c_{02}-c_{03},b_{04},
c_{02}b_{02}^d,c_{02}b_{03}^d\Bigr)
\cap\Bigl(b_{02}^d,c_{03},c_{02}b_{01}\Bigr)\cr
&=\C_1 + (b_{01},b_{02}^d,c_{03}b_{03}^d,
(c_{02}-c_{03})c_{03},b_{04}c_{03}), \cr
}
$$
so that finally $L: b_{03}^d$ equals
$$
\eqalignno{
&\Bigl(b_{02}^d-b_{01}^d,c_{02}-c_{03},b_{01}-b_{04},b_{03}^d-b_{01}^d,
c_{1i}(b_{02}-b_{1i}b_{03}),
c_{1i}c_{1j}(b_{1i}-b_{1j}),
c_{1i}(1-b_{1i}^d)\Bigr)\cr
&\hskip1em\cap
\Bigl(\C_1 + (b_{01},b_{02}^d,c_{03}b_{03}^d,
(c_{02}-c_{03})c_{03},b_{04}c_{03})\Bigr)\cr
&\hskip1em\bigcap_{\Lambda\not=\emptyset}
\Bigl((b_{01},c_{03}b_{04}) +(c_{1i}| i \not \in \Lambda) \cr
&\hskip2em+ (b_{02}-b_{1i}b_{03},b_{1i}-b_{1j},
b_{1i}^d(b_{03}^d,c_{02}-c_{03},b_{04}),
c_{02}b_{1i}^d-c_{03}
| i,j \in \Lambda)\Bigr)\cr
&=
\Bigl(b_{02}^d-b_{01}^d,c_{02}-c_{03},b_{01}-b_{04},b_{03}^d-b_{01}^d,
c_{1i}(b_{02}-b_{1i}b_{03}),
c_{1i}c_{1j}(b_{1i}-b_{1j}),
c_{1i}(1-b_{1i}^d)\Bigr)\cr
&\hskip1em\cap
\Bigl((b_{01},b_{02}^d,c_{03}b_{03}^d,
(c_{02}-c_{03})c_{03},b_{04}c_{03})\cr
&\hskip2em+ (c_{1i}(b_{02}-b_{1i}b_{03}),c_{1i}c_{1j}(b_{1i}-b_{1j}),
c_{1i}b_{1i}^d(c_{02}-c_{03},b_{04}),
c_{1i}(c_{02}b_{1i}^d-c_{03})\Bigr)\cr
&=
(b_{02}^d-b_{01}^d,c_{03}(b_{03}^d-b_{01}^d),
(c_{02}-c_{03})c_{03},b_{04}c_{03}-b_{01}c_{02})\cr
&\hskip2em+ (c_{1i}(b_{02}-b_{1i}b_{03}),c_{1i}c_{1j}(b_{1i}-b_{1j}),
c_{1i}b_{1i}^d(c_{02}-c_{03},b_{04}-b_{01}),
c_{1i}(c_{02}b_{1i}^d-c_{03}))\cr
&\hskip2em+b_{01}(c_{02}-c_{03},b_{01}-b_{04},b_{03}^d-b_{01}^d,
c_{1i}(1-b_{1i}^d)). \cr
}
$$
Thus finally
$$
\eqalignno{
&K : c_{02} b_{03}^d =
\left(c_{01}-c_{02}b_{02}^d,c_{01}-c_{04},s-fb_{02}^d\right)
+ c_{13}(c_{02}-c_{03},b_{01}-b_{04}, b_{03}^d-b_{01}^d)
\cr
&\hskip2em+b_{02}^d (c_{1i}b_{1i}^d - c_{13}b_{13}^d)
+c_{13}b_{03}^d \left( \D_1 + (b_{11}-b_{14},1-b_{1i}^d)\right) \cr
&\hskip2em
+c_{13}b_{03}^d \left(
(c_{11},b_{02},b_{03})(1-b_{2i})
+\left(b_{12}-b_{2i}b_{13},b_{2i}-b_{2j}\right)\right)\delta_{r>2} \cr
&\hskip2em+
\sum_{k=2}^{r-2} b_{03}^d c_{13} \left( \D_k +(1-b_{k+1,i})\right)
+ c_{13}b_{03}^d \left(\D_{r-1} + \C_r\right) \cr
&\hskip2em+
(b_{02}^d-b_{01}^d,c_{03}(b_{03}^d-b_{01}^d),
(c_{02}-c_{03})c_{03},b_{04}c_{03}-b_{01}c_{02})\cr
&\hskip2em+ (c_{1i}(b_{02}-b_{1i}b_{03}),c_{1i}c_{1j}(b_{1i}-b_{1j}),
c_{1i}b_{1i}^d(c_{02}-c_{03},b_{04}-b_{01}),
c_{1i}(c_{02}b_{1i}^d-c_{03}))\cr
&\hskip2em+b_{01}(c_{02}-c_{03},b_{01}-b_{04},b_{03}^d-b_{01}^d,
c_{1i}(1-b_{1i}^d)). \cr
}
$$
By Fact~\factassocses,
$\Ass \left( {R \over K : c_{02}b_{03}^d} \right)
\subseteq
 \Ass \left( {R \over K :c_{02}b_{03}^dc_{13}} \right)
\cup \Ass \left( {R \over (K :c_{02}b_{03}^d) + (c_{13})} \right)$.
Note that
$$
\eqalignno{
&(K : c_{02} b_{03}^d) + (c_{13}) =
(c_{01}-c_{02}b_{02}^d,c_{01}-c_{04},s-fb_{02}^d,c_{13}) \cr
&\hskip2em
+b_{02}^d(c_{1i}b_{1i}^d)
+(b_{02}^d-b_{01}^d,c_{03}(b_{03}^d-b_{01}^d),
(c_{02}-c_{03})c_{03},b_{04}c_{03}-b_{01}c_{02})\cr
&\hskip2em+ (c_{1i}(b_{02}-b_{1i}b_{03}),c_{1i}c_{1j}(b_{1i}-b_{1j}),
c_{1i}b_{1i}^d(c_{02}-c_{03},b_{04}-b_{01}),
c_{1i}(c_{02}b_{1i}^d-c_{03}))\cr
&\hskip2em+b_{01}(c_{02}-c_{03},b_{01}-b_{04},b_{03}^d-b_{01}^d,
c_{1i}(1-b_{1i}^d)). \cr
}
$$
No $b_{13}$, $c_{2i}$ or $b_{2i}$
appear in a minimal generating set of this ideal,
so that by Theorem~\thmnox,
$(K : c_{02} b_{03}^d) + (c_{13})$ gives no new embedded primes of $J$.
Furthermore,
$x$ is a non-zerodivisor on all of these.
Thus it remains to analyze the associated primes of
$K :c_{02}b_{03}^dc_{13}$:
$$
\eqalignno{
&K : c_{02} b_{03}^dc_{13} =
\left(c_{01}-c_{02}b_{02}^d,c_{01}-c_{04},s-fb_{02}^d,
c_{02}-c_{03},b_{01}-b_{04}, b_{03}^d-b_{01}^d\right) \cr
&\hskip2em +b_{03}^d \left( \D_1 + (b_{11}-b_{14},1-b_{1i}^d)\right) \cr
&\hskip2em
+b_{03}^d \left((c_{11},b_{02},b_{03})(1-b_{2i})
+\left(b_{12}-b_{2i}b_{13},b_{2i}-b_{2j}\right)\right)\delta_{r>2} \cr
&\hskip2em+
\sum_{k=2}^{r-2} b_{03}^d \left( \D_k +(1-b_{k+1,i})\right)
+ b_{03}^d \left(\D_{r-1} + \C_r\right) \cr
&\hskip2em+
(b_{02}-b_{13}b_{03},c_{1i}(b_{1i}-b_{13}),
b_{13}^d(c_{02}-c_{03},b_{04}-b_{01}),
c_{02}b_{13}^d-c_{03},b_{01}(1-b_{13}^d))\cr
&\hskip2em+
(b_{02}^d (c_{13}b_{13}^d-c_{1i}b_{1i}^d | i \not = 3) +L' ): c_{13}, \cr
}
$$
where
$$
\eqalignno{
L'&= (b_{02}^d-b_{01}^d,c_{03}(b_{03}^d-b_{01}^d),
(c_{02}-c_{03})c_{03},b_{04}c_{03}-b_{01}c_{02},
c_{1i}c_{1j}(b_{1i}-b_{1j})|i,j \not = 3)\cr
&\hskip2em+ (c_{1i}(b_{02}-b_{1i}b_{03}),
c_{1i}b_{1i}^d(c_{02}-c_{03},b_{04}-b_{01}),
c_{1i}(c_{02}b_{1i}^d-c_{03})|i \not = 3) \cr
&\hskip2em+b_{01}(c_{02}-c_{03},b_{01}-b_{04},b_{03}^d-b_{01}^d,
c_{1i}(1-b_{1i}^d)| i \not = 3).\cr
&= (b_{02}^d-b_{01}^d,c_{03}(b_{03}^d-b_{01}^d),
(c_{02}-c_{03})c_{03},(b_{04}-b_{01})c_{03},
c_{1i}c_{1j}(b_{1i}-b_{1j})|i,j \not = 3)\cr
&\hskip2em+ (c_{1i}(b_{02}-b_{1i}b_{03}),
c_{1i}b_{1i}^d(c_{02}-c_{03},b_{04}-b_{01}),
c_{1i}c_{03}(b_{1i}^d-1)|i \not = 3) \cr
&\hskip2em+b_{01}(c_{02}-c_{03},b_{01}-b_{04},b_{03}^d-b_{01}^d,
c_{1i}(1-b_{1i}^d)| i \not = 3). \cr
}
$$
Clearly
$(b_{02}^d (c_{13}b_{13}^d-c_{1i}b_{1i}^d | i \not = 3) +L' ): c_{13}$ contains
$$
\eqalignno{
&L'' = b_{02}^d (c_{13}b_{13}^d-c_{1i}b_{1i}^d)
+ b_{02}^d b_{13}^d(c_{1i}-c_{1j},b_{02}-b_{1i}b_{03},
1-b_{1i}^d,c_{1i}(b_{1i}-b_{13}) | i,j \not = 3) \cr
&\hskip2em+
(b_{02}^d-b_{01}^d,c_{03}(b_{03}^d-b_{01}^d),
(c_{02}-c_{03})c_{03},(b_{04}-b_{01})c_{03})\cr
&\hskip2em+ (c_{1i}(b_{02}-b_{1i}b_{03}),c_{1i}c_{1j}(b_{1i}-b_{1j}),
c_{1i}b_{1i}^d(c_{02}-c_{03},b_{04}-b_{01}),
c_{1i}c_{03}(b_{1i}^d-1)|i,j \not = 3) \cr
&\hskip2em+b_{01}(c_{02}-c_{03},b_{01}-b_{04},b_{03}^d-b_{01}^d,
c_{1i}(1-b_{1i}^d)| i \not = 3). \cr
}
$$
It turns out that $L'' = L' : c_{13}$,
as the proof below shows.

Let $y \in (b_{02}^d (c_{13}b_{13}^d-c_{1i}b_{1i}^d | i \not = 3) +L' ): c_{13}$.
Write
$$
y c_{13} = \sum_{i\not = 3} y_ib_{02}^d (c_{13}b_{13}^d-c_{1i}b_{1i}^d) + l,
$$
for some $y_i$ in the ring and $l \in L'$.
Then $y_1b_{02}^d c_{11}b_{11}^d \in L' + (c_{12},c_{13},c_{14})$,
so that without loss of generality
$$
y_1 \in (c_{12}, c_{14},b_{02}^d-b_{01}^d,c_{02}-c_{03},b_{01}-b_{04},
b_{03}^d-b_{01}^d,
b_{02}-b_{11}b_{03},
b_{11}^d-1).
$$
Thus $y_1 b_{02}^d (c_{13}b_{13}^d-c_{11}b_{11}^d)$ is contained in
$$
\eqalignno{
&b_{02}^d (c_{13}b_{13}^d-c_{11}b_{11}^d)
\left(c_{1j},c_{02}-c_{03},b_{01}-b_{04},b_{03}^d-b_{01}^d,
b_{02}-b_{11}b_{03},b_{11}^d-1,| j \not = 1,3\right) \cr
&\subseteq  L' + b_{02}^d (c_{13}b_{13}^d-c_{11}b_{11}^d)
\left(c_{1j}, b_{02}-b_{11}b_{03},b_{11}^d-1| j \not = 1,3\right) \cr
&\subseteq  L' +
b_{02}^d (c_{13}b_{13}^d-c_{11}b_{11}^d) (c_{1j} | j \not = 1,3)
+b_{02}^d c_{13}b_{13}^d(b_{02}-b_{11}b_{03},b_{11}^d-1) \cr
&\subseteq  L' +
b_{02}^d c_{11}(c_{13}b_{13}^d-c_{1j}b_{1j}^d)) + b_{02}^d c_{13}b_{13}^d (c_{1j}-c_{11},
b_{02}-b_{11}b_{03},b_{11}^d-1| j \not = 1,3). \cr
}
$$
Thus for some $y' \in b_{02}^d c_{13}b_{13}^d (c_{1j}-c_{11},
b_{02}-b_{11}b_{03},b_{11}^d-1| j \not = 1,3) \subseteq L''$
and some $y'_2, y'_4$ in the ring,
$$
(y-y')c_{13} - \sum_{i=2,4} y'_i b_{02}^d (c_{13}b_{13}^d-c_{1i}b_{1i}^d) \in L'.
$$
Then $y'_2b_{02}^d c_{12}b_{12}^d$
is in $L' + (c_{13},c_{14})$,
so that
$$
\eqalignno{
y'_2& c_{12}b_{12}^d
\in
(c_{13},c_{14},b_{02}^d-b_{01}^d,c_{02}-c_{03},b_{01}-b_{04},
b_{03}^d-b_{01}^d) \cr
&\hskip2em+ (c_{1i}(1-b_{1i}^d),
c_{1i}(b_{02}-b_{1i}b_{03}),c_{11}c_{12}(b_{11}-b_{12})
| i=1,2), \cr
}
$$
whence
$$
y'_2\in
(c_{13},c_{14},c_{02}-c_{03},b_{01}-b_{04},b_{03}^d-b_{01}^d,
1-b_{12}^d,b_{02}-b_{12}b_{03},
c_{11}(b_{11}-b_{12})).
$$
By similar reasoning as for $y_1$,
there exists $y'' \in L"$ and $y''_4$ in the ring
such that
$(y-y'-y'')c_{13} -y''_4b_{02}^d (c_{13}b_{13}^d-c_{14}b_{14}^d) \in L'$.
Then $y''_4b_{02}^d c_{14}b_{14}^d \in L' + (c_{13})$,
and again one can conclude that
$y''_4\in L''+ (c_{13})$.
Thus $y$ is an element of $L'$ modulo $L''$,
so that $L'' = L' : c_{13}$.
Thus finally
$$
\eqalignno{
&K : c_{02} b_{03}^dc_{13} =
\left(c_{01}-c_{02}b_{02}^d,c_{01}-c_{04},s-fb_{02}^d,
c_{02}-c_{03},b_{01}-b_{04}, b_{03}^d-b_{01}^d\right) \cr
&\hskip2em +b_{03}^d \left( \D_1 + (b_{11}-b_{14},1-b_{1i}^d)\right) \cr
&\hskip2em
+b_{03}^d \left((c_{11},b_{02},b_{03})(1-b_{2i})
+\left(b_{12}-b_{2i}b_{13},b_{2i}-b_{2j}\right)\right)\delta_{r>2} \cr
&\hskip2em+
\sum_{k=2}^{r-2} b_{03}^d \left( \D_k +(1-b_{k+1,i})\right)
+ b_{03}^d \left(\D_{r-1} + \C_r\right) \cr
&\hskip2em+
(b_{02}-b_{13}b_{03},c_{1i}(b_{1i}-b_{13}),
b_{02}^d b_{13}^d(b_{1i}-b_{13})b_{03},
b_{02}^d-b_{01}^d)\cr
&\hskip2em+ (c_{02},b_{01})(1-b_{13}^d,c_{1i}c_{02}(1-b_{1i}^d)|i \not = 3). \cr
}
$$
By Fact~\factassocses,
$\Ass \left( {R \over K : c_{02}b_{03}^dc_{13}} \right)
\subseteq
 \Ass \left( {R \over K :c_{02}b_{03}^{2d}c_{13}} \right)
\cup \Ass \left( {R \over (K :c_{02}b_{03}^dc_{13}) + (b_{03}^d)} \right)$.
The latter ideal equals and decomposes as:
$$
\eqalignno{
&(K :c_{02}b_{03}^dc_{13}) + (b_{03}^d) =
(s,c_{01},c_{04},
c_{02}-c_{03},b_{01}-b_{04},b_{01}^d,b_{02}-b_{13}b_{03},
b_{03}^d) \cr
&\hskip2em+
(c_{1i}(b_{1i}-b_{13}),
c_{02}(b_{13}^d-1),b_{01}(1-b_{13}^d),
c_{02}c_{1i}(b_{1i}^d-1)),b_{01}(c_{1i}(1-b_{1i}^d)) \cr
&=\Bigl(
\C_0 + (s,b_{01},b_{04},b_{02}-b_{13}b_{03},
b_{03}^d,c_{1i}(b_{1i}-b_{13})) \Bigr) \cr
&\hskip1em\cap\Bigl(
(s,c_{01},c_{04},
c_{02}-c_{03},b_{01}-b_{04},b_{01}^d,b_{02}-b_{13}b_{03},
b_{03}^d) \cr
&\hskip2em+
(c_{1i}(b_{1i}-b_{13}),b_{13}^d-1,
c_{1i}(b_{1i}^d-1)) \Bigr), \cr
}
$$
which is an intersection of $Q_{3\Lambda}$-
and $Q_{2\Lambda\alpha}$-primary components,
where $\Lambda$ varies over all the subsets of $\{1,2,3,4\}$
for which $3 \in \Lambda$.
These do not give any new embedded primes of $J$,
and furthermore none of these primes contains $x$.

It remains to analyze the associated primes of $K :c_{02}b_{03}^{2d}c_{13}$:
$$
\eqalignno{
&K : c_{02} b_{03}^{2d}c_{13} =
\left(c_{01}-c_{02}b_{02}^d,c_{01}-c_{04},s-fb_{02}^d,
c_{02}-c_{03},b_{01}-b_{04}, b_{02}-b_{13}b_{03}\right) \cr
&\hskip2em + \D_1 + (b_{11}-b_{14},1-b_{1i}^d)
+ \left((c_{11},b_{02},b_{03})(1-b_{2i})
+\left(b_{12}-b_{2i}b_{13},b_{2i}-b_{2j}\right)\right)\delta_{r>2} \cr
&\hskip2em+
\sum_{k=2}^{r-2} \left( \D_k +(1-b_{k+1,i})\right) + \D_{r-1} + \C_r
+ (b_{13}^d(b_{1i}-b_{13})b_{03}) \cr
&\hskip2em+
\left((b_{03}^d-b_{01}^d,c_{1i}(b_{1i}-b_{13}))
+ (c_{03},b_{01})(c_{1i}(b_{1i}^d-1) ,1-b_{13}^d)\right) : b_{03}^d \cr
&=
\left(c_{01}-c_{02}b_{02}^d,c_{01}-c_{04},s-fb_{02}^d,
c_{02}-c_{03},b_{01}-b_{04}, b_{02}-b_{13}b_{03}\right) \cr
&\hskip2em + \D_1 + (b_{11}-b_{14},1-b_{1i}^d)
+ \left((c_{11},b_{02},b_{03})(1-b_{2i})
+\left(b_{12}-b_{2i}b_{13},b_{2i}-b_{2j}\right)\right)\delta_{r>2} \cr
&\hskip2em+
\sum_{k=2}^{r-2} \left( \D_k +(1-b_{k+1,i})\right) + \D_{r-1} + \C_r+
(b_{1i}-b_{13})(b_{03},c_{1i})+ \left(b_{03}^d-b_{01}^d\right).\cr
}
$$
Note that this decomposes as
$$
\eqalignno{
&\Bigl((c_{01}-c_{02}b_{02}^d,c_{01}-c_{04},s-fb_{02}^d,
c_{02}-c_{03},b_{01}-b_{04},b_{02}-b_{13}b_{03}) + \D_1\cr
&\hskip1em
+ (1-b_{1i}^d, (1-b_{2i})\delta_{r > 2},b_{13}-b_{1i},b_{03}^d-b_{01}^d)
+ \sum_{k=2}^{r-2}
\left( \D_k +(1-b_{k+1,i})\right) +\D_{r-1} + \C_r \Bigr) \cr
&\cap\Bigl((s,c_{01},c_{04},
c_{02}-c_{03},b_{01}-b_{04},b_{02},b_{03},b_{01}^d)+\C_1
+ (b_{11}-b_{14},1-b_{1i}^d) \cr
&\hskip1em+(b_{12}-b_{2i}b_{13},b_{2i}-b_{2j}) \delta_{r > 2}
+\sum_{k=2}^{r-2}
\left( \D_k +(1-b_{k+1,i})\right) +\D_{r-1} + \C_r \Bigr) \cr
}
$$
(The key to this decomposition is the fact that
$(c_{11},b_{03})$ intersected with the first component
is contained in $K : c_{02}b_{03}^{2d}c_{13}$.)
The first component above is $p_r$,
so all of its associated primes are minimal over $J$.
It is easy to read off the associated primes of the last component as well.
First note that none of these primes contain $x$,
which finishes the proof that $x$ is  a non-zerodivisor modulo $K$.
Thus as $J \subseteq K$ and $xK \subseteq J$,
it follows that $K$ equals $J : x$.

It remains to determine the associated prime ideals of this last
component of $K : c_{02}b_{03}^{2d}c_{13}$ in the display above.
The last component is the intersection
of $Q_{4r\alpha\beta\gamma}$-primary components,
as $\alpha, \beta$, and $\gamma$ vary over all the $d$th roots of unity in $K$.
Note that $Q_{4,r\alpha\alpha\alpha}$-component contains $p_r$
and is thus redundant in the decomposition.
But coloning with $b_{13}-b_{1i}$ for various $i$
shows that the remaining prime ideals
are indeed associated to $K$ and thus to $J$.

This proves

\thm
\label{\thmthirdemb}
Let $n > 2$.
For $r \in \{2, \ldots, n-2\}$,
set $x = f^3(c_{21}\cdots c_{r-1,1})$
$(b_{13}\cdots b_{r-1,3})$
$c_{r+1,1}(1-b_{r1})$,
and for $r = n-1, n$,
set $x = f^3(c_{21}\cdots c_{r-1,1})(b_{13}\cdots b_{r-1,3})$.
Then the set of embedded primes of $J$ not containing $x$ is contained in
$$
\eqalignno{
&\lbrace Q_{1\Lambda}, Q_{2\Lambda'\alpha}, Q_{3\Lambda'} |
\Lambda, \Lambda' \subseteq \{1,2,3,4\}, |\Lambda'| > 0 \rbrace \cr
&\hskip2em
\cup \lbrace
Q_{4,2\alpha\beta}\delta_{n=2},
Q_{4r\alpha\beta\gamma} \delta_{n>2} |
r = 2, \ldots, n; \alpha^d =\beta^d =\gamma^d = 1,
|\{\alpha, \beta, \gamma\}| > 1 \rbrace,
}
$$
and each listed prime ideal is associated to $J$.
\endb

These new associated primes are also recorded in a table:

$$
\displayindent=1em
\halign{\mvrule \hskip 0.1em #\hfil \hskip 0.1em
& \mvrule \hskip 0.5em \relax #\hfil \mvrule \cr
\noalign{\hrule}
embedded prime ($\alpha^d = \beta^d = \gamma^d = 1$,
$\alpha, \beta, \gamma$ not all equal) & height \cr
\ \ \ $n > 2$, $r = 2, \ldots, n$ & \cr
\noalign{\hrule}
\noalign{\hrule}
$Q_{4r\alpha\beta\gamma} =
\left(s, c_{01}, c_{03}-c_{02}, c_{04}, b_{01},b_{02},b_{03},b_{04}\right)$
& $7r+2+4\delta_{r < n}$ \cr
\hskip1em$+ (b_{12}-b_{2i}b_{13},b_{2i}-b_{2j}) \delta_{r > 2}
+(b_{11}-\alpha,b_{14}-\alpha,b_{12}-\beta,b_{13}-\gamma)$
& \cr
\hskip1em$+\C_1 +\D_2 + \cdots + \D_{r-1} + \C_r + B_{3,r-1}$
& \cr
\noalign{\hrule}
}
$$

\section{Reduction to $\bf(J(n,d) : sc_{02}) + (c_{02},f)$}
\sectlabel{\sectreduction}
\pagelabel{\qsectreduction}

In this section the finding of the embedded primes of $J$
gets reduced to that of finding the associated primes of ideas
on which recursion can be applied.
The main methods
are again repeated applications of Facts~\factdec\ and \factassocses.
For example,
the set of associated primes of $J$
is contained in $\Ass \left({R \over J + (s)}\right)
\cup \Ass\left({R \over J : s}  \right)$.

To start off,
the decomposition of $J + (s)$ is easy:
$$
\eqalignno{
J + (s) &= (s) + f\left(c_{01},c_{04},c_{02}b_{02}^d,c_{03}b_{03}^d,
c_{02}b_{01}-c_{03}b_{04},
c_{02}c_{1i}(b_{02}-b_{1i}b_{03})\right) \cr
&= (s,f) \cap
\left(s,c_{01},c_{04},c_{02}b_{02}^d,c_{03}b_{03}^d,
c_{02}b_{01}-c_{03}b_{04},
c_{02}c_{1i}(b_{02}-b_{1i}b_{03})\right) \cr
&= p_{-1} \cap
\left(s,c_{01},c_{04},c_{02}b_{02}^d,c_{03}b_{03}^d,
c_{02}b_{01}-c_{03}b_{04},
c_{02}c_{1i}(b_{02}-b_{1i}b_{03}),c_{02}\right) \cr
&\hskip 3em \cap
\left(\left(s,c_{01},c_{04},c_{02}b_{02}^d,c_{03}b_{03}^d,
c_{02}b_{01}-c_{03}b_{04},
c_{02}c_{1i}(b_{02}-b_{1i}b_{03})\right) : c_{02}\right) \cr
&= p_{-1} \cap
\left(s,c_{01},c_{04},c_{02},c_{03}b_{03}^d,c_{03}b_{04} \right) \cr
&\hskip 3em \cap
\left(s,c_{01},c_{04},b_{02}^d,c_{03}b_{03}^d,
c_{02}b_{01}-c_{03}b_{04},b_{01}b_{03}^d,
c_{1i}(b_{02}-b_{1i}b_{03})\right) \cr
&= p_{-1} \cap
\left(s,c_{01},c_{04},c_{02},c_{03} \right) \cap
\left(s,c_{01},c_{04},c_{02},b_{03}^d,b_{04} \right) \cr
&\hskip 3em \cap
\left(\left(s,c_{01},c_{04},b_{02}^d,c_{03}b_{03}^d,
c_{02}b_{01}-c_{03}b_{04},b_{01}b_{03}^d,
c_{1i}(b_{02}-b_{1i}b_{03})\right) : b_{03}^d\right) \cr
&\hskip 3em \cap
\left(s,c_{01},c_{04},b_{02}^d,c_{03}b_{03}^d,
c_{02}b_{01}-c_{03}b_{04},b_{01}b_{03}^d,
c_{1i}(b_{02}-b_{1i}b_{03}),b_{03}^d\right) \cr
&= p_{-1} \cap \left(s,c_{01},c_{04},c_{02},c_{03} \right) \cap p_{-2} \cr
&\hskip 3em \cap
\left(s,c_{01},c_{04},b_{02}^d,c_{03},b_{01},
c_{1i}(b_{02}-b_{1i}b_{03}),
c_{1i}c_{1j}(b_{1i}-b_{1j}),c_{1i}b_{1i}^d \right) \cr
&\hskip 3em \cap
\left(s,c_{01},c_{04},b_{02}^d,b_{03}^d,
c_{02}b_{01}-c_{03}b_{04},
c_{1i}(b_{02}-b_{1i}b_{03}) \right) \cr
&= p_{-1} \cap \left(s,c_{01},c_{04},c_{02},c_{03}\right)
\cap p_{-2}\cap p_{-4}\cr
&\hskip 3em \cap
\left(s,c_{01},c_{04},b_{02}^d,b_{03}^d,
c_{02}b_{01}-c_{03}b_{04},
c_{1i}(b_{02}-b_{1i}b_{03}),c_{1i}c_{1j}(b_{1i}-b_{1j}) \right) \cr
&\hskip 3em \cap
\left(s,c_{01},c_{04},b_{02},b_{03},
c_{02}b_{01}-c_{03}b_{04} \right) \cr
&= p_{-1} \cap
\left(s,c_{01},c_{04},c_{02},c_{03} \right) \cap
p_{-2} \cap p_{-4} \cap q_1 \cap p_{-3}. \cr
}
$$
Recall that $p_{-1}, p_{-2}$ and $p_{-3}$ are minimal components of $J$,
that $p_{-4}$ and $q_1$ are the intersections of 16 components of $J$ each,
but that by Theorem~\thmnox,
$\left(s,c_{01},c_{04},c_{02},c_{03} \right)$ is not associated to $J$
as it is not minimal.

This proves (by Fact~\factassocses):

\thm
The set of embedded primes of $J$ equals
$\lbrace Q_{1\Lambda} | \Lambda \rbrace \cup
\Ass \left( {R \over J : s} \right)$.~\qed
\endb

The next task is to compute $J : s$
and to analyze its associated primes.
Any associated prime of $J:s$ is also associated to $J$.
Computing $J:s$ is straightforward (next Theorem),
but analyzing its associated primes takes many steps
and the rest of this paper.

\def\JJ{J_2}

\thm
\label{\thmcolons}
Let $\JJ$ be the ideal in $R$ generated by all the $h_{rj}/s$, $r \ge 2$.
(Note that all these $h_{rj}$ are multiples of $s$.)
Then
$J:s$ equals
$$
\eqalignno{
J : s &=
\left(
c_{01}-c_{02}b_{01}^d, c_{03}-c_{02}, c_{04}-c_{02}b_{04}^d\right) + \JJ
+ c_{02}(fb_{01}^d-s) \cr
&\hskip2em +(f c_{02}, c_{02}^2)\left(
b_{01}^d-b_{02}^d, b_{04}^d-b_{03}^d,
b_{01}-b_{04}\right) \cr
&\hskip2em
+c_{02} \left(
b_{01}b_{03}^d-b_{04}b_{02}^d,
c_{1i}(b_{02}-b_{1i}b_{03})\right) \cr
&\hskip2em+ c_{02}\left(
c_{1i}c_{1j}(b_{1i}-b_{1j}),
c_{1i}(b_{01}-b_{1i}^db_{04}),
c_{02}c_{1i}(1-b_{1i}^d) \right). \cr
}
$$
where the indices $i$ and $j$ vary from $1$ to $4$.
\endb

\proof
First observe that
$$
J = s \left(
c_{01}-c_{02}b_{01}^d, c_{03}-c_{02}, c_{04}-c_{02}b_{04}^d\right) + s\JJ +
fK + (fc_{01}-sc_{02}),
$$
where
$K =
\left(
c_{01}-c_{02}b_{02}^d, c_{04}-c_{03}b_{03}^d,
c_{01}-c_{04},
c_{02}b_{01}-c_{03}b_{04},
c_{02}c_{1i}(b_{02}-b_{1i}b_{03})\right)$.
Thus $J : s = \left(
c_{01}-c_{02}b_{01}^d, c_{03}-c_{02}, c_{04}-c_{02}b_{04}^d\right) + \JJ$
$+ (fK + (fc_{01}-sc_{02})): s$.
Let $x \in (fK + (fc_{01}-sc_{02})) : s$.
Write $xs = kf + a(fc_{01}-sc_{02})$
for some $k \in K$ and $a \in R$.
By adding to $x$ a multiple of $fc_{01}-sc_{02}$ and an element of $fK$,
and correspondingly changing $a$ and $k$,
without loss of generality no $s$ appears in $a$,
and as $fK \cap (s) = sfK$,
without loss of generality also no $s$ appears in $k$.
 From $xs = kf + a(fc_{01}-sc_{02})$ it follows that
$$
a \in (K + (s)) : fc_{01} = (s) + (K : c_{01}),
$$
and as no $s$ appears in $a$ and the generators of $K$,
actually $a \in K : c_{01}$.
By Fact~\factfewvars,
$K : c_{01} = K : c_{02}b_{02}^d=
(c_{01}-c_{02}b_{02}^d, c_{04}-c_{03}b_{03}^d) + (K' : c_{02}b_{02}^d)$,
where
$$
K' = \left( c_{02}b_{02}^d-c_{03}b_{03}^d,
c_{02}b_{01}-c_{03}b_{04}, c_{02}c_{1i}(b_{02}-b_{1i}b_{03})\right).
$$
Then
$$
\eqalignno{
K' : c_{02} &=
\left( c_{02}b_{02}^d-c_{03}b_{03}^d, c_{02}b_{01}-c_{03}b_{04}\right) : c_{02}
+ \left(c_{1i}(b_{02}-b_{1i}b_{03})\right) \cr
&= \left( c_{02}b_{02}^d-c_{03}b_{03}^d, c_{02}b_{01}-c_{03}b_{04},
b_{01}b_{03}^d-b_{04}b_{02}^d, c_{1i}(b_{02}-b_{1i}b_{03})\right), \cr
}
$$
and by the same proof as on page~\qcalcc,
$$
\eqalignno{
K' : c_{02}b_{02}^d &=
\left( c_{02}b_{02}^d-c_{03}b_{03}^d, c_{02}b_{01}-c_{03}b_{04},
b_{01}b_{03}^d-b_{04}b_{02}^d, c_{1i}(b_{02}-b_{1i}b_{03})\right) \cr
&\hskip2em +\left(c_{1i}c_{1j}(b_{1i}-b_{1j}),
c_{1i}(b_{01}-b_{1i}^db_{04}), c_{1i}(c_{03}-b_{1i}^dc_{02}) \right). \cr
}
$$
Thus
$$
\eqalignno{
a \in K : c_{01} &=
\left(c_{01}-c_{02}b_{02}^d, c_{04}-c_{03}b_{03}^d\right) \cr
&+ \left( c_{02}b_{02}^d-c_{03}b_{03}^d,
c_{02}b_{01}-c_{03}b_{04},
b_{01}b_{03}^d-b_{04}b_{02}^d,
c_{1i}(b_{02}-b_{1i}b_{03})\right) \cr
&+ \left(
c_{1i}c_{1j}(b_{1i}-b_{1j}),
c_{1i}(b_{01}-b_{1i}^db_{04}),
c_{1i}(c_{03}-b_{1i}^dc_{02}) \right). \cr
}
$$

Recall that $x \in K : s$
and $sx = kf + a(fc_{01}-sc_{02})$
for some $k \in K$ and $a \in K: c_{01}$.
Thus $s(x+ac_{02}) = f(k+ac_{01})$,
and as no $s$ appears in $a$ and in $k$,
$x+ac_{02} = 0$,
so that $x \in c_{02} (K:c_{01})$.
Thus
$$
\eqalignno{
J : s &\subseteq
\left(
c_{01}-c_{02}b_{01}^d, c_{03}-c_{02}, c_{04}-c_{02}b_{04}^d\right) + \JJ
 + (fc_{01}-sc_{02}) + fK \cr
&\hskip2em
+c_{02}
\left(c_{01}-c_{02}b_{02}^d, c_{04}-c_{03}b_{03}^d\right) \cr
&\hskip2em
+c_{02} \left( c_{02}b_{02}^d-c_{03}b_{03}^d,
c_{02}b_{01}-c_{03}b_{04},
b_{01}b_{03}^d-b_{04}b_{02}^d,
c_{1i}(b_{02}-b_{1i}b_{03})\right) \cr
&\hskip2em+ c_{02}\left(
c_{1i}c_{1j}(b_{1i}-b_{1j}),
c_{1i}(b_{01}-b_{1i}^db_{04}),
c_{1i}(c_{03}-b_{1i}^dc_{02}) \right) \cr
&=
\left(
c_{01}-c_{02}b_{01}^d, c_{03}-c_{02}, c_{04}-c_{02}b_{04}^d\right) + \JJ
+ c_{02}(fb_{01}^d-s) \cr
&\hskip2em +f \left(
c_{01}-c_{02}b_{02}^d, c_{04}-c_{03}b_{03}^d,
c_{01}-c_{04},
c_{02}b_{01}-c_{03}b_{04},
c_{02}c_{1i}(b_{02}-b_{1i}b_{03})\right) \cr
&\hskip2em
+c_{02}^2
\left(b_{01}^d-b_{02}^d, b_{04}^d-b_{03}^d,b_{02}^d-b_{03}^d,b_{01}-b_{04},
c_{1i}(1-b_{1i}^d) \right) \cr
&\hskip2em
+c_{02} \left(
b_{01}b_{03}^d-b_{04}b_{02}^d,
c_{1i}(b_{02}-b_{1i}b_{03}),
c_{1i}c_{1j}(b_{1i}-b_{1j}),
c_{1i}(b_{01}-b_{1i}^db_{04})\right) \cr
&=
\left(
c_{01}-c_{02}b_{01}^d, c_{03}-c_{02}, c_{04}-c_{02}b_{04}^d\right) + \JJ
+ c_{02}(fb_{01}^d-s) \cr
&\hskip2em +f c_{02}\left(
b_{01}^d-b_{02}^d, b_{04}^d-b_{03}^d,
b_{01}-b_{04}\right) \cr
&\hskip2em
+c_{02}^2
\left(b_{01}^d-b_{02}^d, b_{04}^d-b_{03}^d,b_{01}-b_{04},
c_{1i}(1-b_{1i}^d) \right) \cr
&\hskip2em
+c_{02} \left(
b_{01}b_{03}^d-b_{04}b_{02}^d,
c_{1i}(b_{02}-b_{1i}b_{03}),
c_{1i}c_{1j}(b_{1i}-b_{1j}),
c_{1i}(b_{01}-b_{1i}^db_{04}) \right). \cr
}
$$
It is easy to verify that the other inclusion also holds,
which proves the theorem.
\qed

Incidentally,
this also shows:

\prop
\label{\proprad}
The Mayr-Meyer ideal $J(n,d)$ is not a radical ideal:
the element $sc_{02}(b_{01}-b_{04})$ is in $\sqrt J$ but not in $J$.
\qed
\endb

This was already proved in [S2] with the assumption that $d \ge 2$,
without giving an element of the radical which is not in the ideal.

Furthermore,
it is easy to see the following:

\cor
Let $a$ be one of the listed generators of $J(n,d): s$.
Then $s \cdot a$ can be written as a linear combination of the generators of $J(n,d)$
with coefficients of degree at most $2d+1$.
Also,
$c_{02} b_{01}^d c_{11} \cdots c_{n-2,1} (c_{n-1,1} - c_{n-1,4})$
lies in $J_(n,d): s$. 
\qed
\endb

Let $\JJ'$ be the ideal obtained from $\JJ$
after rewriting each
$c_{01}$ as $c_{02}b_{01}^d$,
$c_{03}$ as $c_{02}$,
and $c_{04}$ as $c_{02}b_{04}^d$.
Note that
$\JJ'$ is a multiple of $c_{02}$
and that the theorem above also holds with $\JJ'$ in place of $\JJ$.

Observe that $(J : s) + (c_{02}) = \C_0 = p_0$,
a minimal prime ideal over $J$.
Thus by Fact~\factassocses:

\thm
\label{\thmfirstemb}
The set of embedded primes of $J$ equals
$\lbrace Q_{1\Lambda} | \Lambda \rbrace \cup
\Ass \left( {R \over J : sc_{02}} \right)$,
and that in turn is contained in
$$
\lbrace Q_{1\Lambda} | \Lambda \rbrace
\cup \Ass \left( {R \over J : sc_{02}^2} \right)
\cup \Ass \left( {R \over (J : sc_{02}) + (c_{02})} \right),
$$
where
$$
\Ass \left( {R \over (J : sc_{02}) + (c_{02})} \right)
\subseteq
\Ass \left( {R \over ((J : sc_{02}) + (c_{02})):f} \right)
\cup \Ass \left( {R \over (J : sc_{02}) + (c_{02},f)} \right). 
$$
\endb

Here are all the ideals appearing in this theorem:
$$
\eqalignno{
&J : sc_{02} =
\left(
c_{01}-c_{02}b_{01}^d, c_{03}-c_{02}, c_{04}-c_{02}b_{04}^d\right) + \JJ'/c_{02}
+ (fb_{01}^d-s) \cr
&\hskip2em +(f, c_{02})\left(
b_{01}^d-b_{02}^d, b_{04}^d-b_{03}^d, b_{01}-b_{04}\right)
+ \left( b_{01}b_{03}^d-b_{04}b_{02}^d \right) \cr
&\hskip2em+ \left(
c_{1i}(b_{02}-b_{1i}b_{03}),
c_{1i}c_{1j}(b_{1i}-b_{1j}),
c_{1i}(b_{01}-b_{1i}^db_{04}),
c_{02}c_{1i}(1-b_{1i}^d) \right), \cr
&J : sc_{02}^2 =
\left(
c_{01}-c_{02}b_{01}^d, c_{03}-c_{02}, c_{04}-c_{02}b_{04}^d\right) + \JJ'/c_{02}
+ (fb_{01}^d-s) \cr
&\hskip2em +\left(
b_{01}^d-b_{02}^d, b_{04}^d-b_{03}^d, b_{01}-b_{04},
c_{1i}(b_{02}-b_{1i}b_{03}), c_{1i}c_{1j}(b_{1i}-b_{1j}),
c_{1i}(1-b_{1i}^d) \right), \cr
&(J : sc_{02}) + (c_{02}) = \C_0 + \JJ'/c_{02}
+ (fb_{01}^d-s) +f\left(
b_{01}^d-b_{02}^d, b_{04}^d-b_{03}^d, b_{01}-b_{04}\right) \cr
&\hskip2em
+ \left(
b_{01}b_{03}^d-b_{04}b_{02}^d,
c_{1i}(b_{02}-b_{1i}b_{03}),
c_{1i}c_{1j}(b_{1i}-b_{1j}),
c_{1i}(b_{01}-b_{1i}^db_{04}) \right), \cr
&((J : sc_{02}) + (c_{02})):f = \C_0 + \JJ'/c_{02}
+ \left(fb_{01}^d-s,
b_{01}^d-b_{02}^d, b_{04}^d-b_{03}^d, b_{01}-b_{04}\right)\cr
&\hskip2em
+ \left(
c_{1i}(b_{02}-b_{1i}b_{03}), c_{1i}c_{1j}(b_{1i}-b_{1j}),
c_{1i}(b_{01}-b_{1i}^db_{04}) \right), \cr
&(J : sc_{02}) + (c_{02},f) = \C_0 + \JJ'/c_{02}+
\left(s,f, b_{01}b_{03}^d-b_{04}b_{02}^d\right) \cr
&\hskip2em+ \left(
c_{1i}(b_{02}-b_{1i}b_{03}),
c_{1i}c_{1j}(b_{1i}-b_{1j}),
c_{1i}(b_{01}-b_{1i}^db_{04}) \right). \cr
}
$$
\pagelabel{\qinought}
Observe that $(J : sc_{02}^2) + (b_{01}^d)$ equals
$$
\left(
c_{01}, c_{03}-c_{02}, c_{04}, s, b_{0i}^d, b_{01}-b_{04},
c_{1i}(b_{02}-b_{1i}b_{03}), c_{1i}c_{1j}(b_{1i}-b_{1j}),
c_{1i}(1-b_{1i}^d) \right),
$$
and $(((J : sc_{02}) + (c_{02})):f) + (b_{01}^d)$ equals
$$
\C_0
+ \left( s, b_{0i}^d, b_{01}-b_{04}\right)
+ \left(
c_{1i}(b_{02}-b_{1i}b_{03}), c_{1i}c_{1j}(b_{1i}-b_{1j}),
c_{1i}(b_{01}-b_{1i}^db_{04}) \right).
$$
Clearly the associated primes of these two ideals do not contain $x$,
where $x =fc_{21}b_{13}(b_{21}-b_{22})$ when $n > 2$
and $x = f$ otherwise.
Thus by Theorem~\thmnox,
these ideals do not contribute anything new to the set of embedded primes of $J$.

Thus by another application of Fact~\factassocses,

\thm
\label{\thmsecondemb}
The set of embedded primes of $J$ is contained in
$$
\eqalignno{
&\lbrace Q_{1\Lambda}, Q_{2\Lambda'\alpha}, Q_{3\Lambda'} |
\Lambda, \Lambda' \subseteq \{1,2,3,4\}, |\Lambda'| > 0, \alpha^d = 1 \rbrace
\cr
&\cup \Ass \left( {R \over J : sc_{02}^2b_{01}^d} \right)
\cup \Ass \left( {R \over ((J : sc_{02})\hskip-0.1em+\hskip-0.1em (c_{02})):fb_{01}^d} \right)
\cup \Ass \left( {R \over (J : sc_{02})\hskip-0.1em+\hskip-0.1em(c_{02},f)} \right).
\eqed \cr
}
$$
\endb

The rest of this section determines the associated prime ideals of
$J : sc_{02}^2b_{01}^d$ and $((J : sc_{02})+(c_{02})):fb_{01}^d$.

Define $\JJ''$ to be the ideal
$$
\JJ'' = (h_{rj} | r \ge 2) \hbox{\ \ with setting $s = c_{01} = c_{04} = 1$}.
$$
This is the same as taking the ideal $\JJ$,
rewriting each $c_{01}$ and $c_{04}$ as $c_{02}b_{01}^d$
(whence each element is divisible by $c_{02}b_{01}^d$),
and then dividing that ideal by $c_{02}b_{01}^d$.
Recall that $\JJ'$ is the ideal obtained from $J_2$
by rewriting each $c_{01}$ as $c_{02}b_{01}^d$ and
$c_{04}$ as $c_{02}b_{04}^d$.
Then
$$
\JJ'' = \D_1 +
\sum_{r=1}^{n-1} c_{11} \cdots c_{r1}
\Bigl( \D_{r+1} + (b_{r1}-b_{r4}, c_{r+1,1}(b_{r2}-b_{r+1,i}b_{r3})) \Bigr),
$$
using the convention that $c_{ni}=1 = b_{ni}$ and $\D_n = (0)$,
and
$$
\JJ'/c_{02} + \left(b_{01}-b_{04}\right) =
\JJ''b_{01}^d + \left(b_{01}-b_{04}\right).
$$
Thus
$$
\eqalignno{
&J : sc_{02}^2 b_{01}^d = \JJ'' +
\left(
c_{01}-c_{02}b_{01}^d, c_{03}-c_{02}, c_{04}-c_{02}b_{04}^d,s-fb_{01}^d\right)
\cr
&\hskip2em +\left(
b_{01}^d-b_{02}^d, b_{04}^d-b_{03}^d, b_{01}-b_{04},
c_{1i}(b_{02}-b_{1i}b_{03}), c_{1i}c_{1j}(b_{1i}-b_{1j}),
c_{1i}(1-b_{1i}^d) \right): b_{01}^d \cr
&\hskip1em= \JJ'' +
\left(
c_{01}-c_{02}b_{01}^d, c_{03}-c_{02}, c_{04}-c_{02}b_{04}^d,s-fb_{01}^d\right)
\cr
&\hskip2em +\left(
b_{01}^d-b_{02}^d, b_{04}^d-b_{03}^d, b_{01}-b_{04}\right)
+ c_{11}\left( b_{02}-b_{1i}b_{03}, c_{11}(b_{1i}-b_{1j}), 1-b_{1i}^d \right), \cr
}
$$
and
$$
\eqalignno{
&((J : sc_{02}) + (c_{02})):fb_{01}^d =
 \C_0 + \JJ''+ (s-fb_{01}^d, b_{01}-b_{04}) \cr
&\hskip2em+\left(
b_{01}^d-b_{02}^d, b_{01}^d-b_{03}^d,
c_{1i}(b_{02}-b_{1i}b_{03}), c_{1i}c_{1j}(b_{1i}-b_{1j}),
c_{1i}b_{01}(1-b_{1i}^d) \right): b_{01}^d \cr
&\hskip1em=
\C_0 + \JJ''+
\left(s-fb_{01}^d,b_{01}^d-b_{02}^d, b_{04}^d-b_{03}^d, b_{01}-b_{04}\right)\cr
&\hskip2em
+ c_{11} \left(b_{02}-b_{1i}b_{03}, c_{11}(b_{1i}-b_{1j}),1-b_{1i}^d \right). \cr
}
$$
Let $L$ be either of the two ideals.
Then $L$ is of the form
$$
L_0 +\JJ'' +
(s-fb_{01}^d,b_{01}^d-b_{02}^d, b_{04}^d-b_{03}^d, b_{01}-b_{04})
+ c_{11}(b_{02}-b_{1i}b_{03}, c_{11}(b_{1i}-b_{1j}),1-b_{1i}^d),
$$
where $L_0$ is either $\C_0$ or
$(c_{01}-c_{02}b_{01}^d, c_{03}-c_{02}, c_{04}-c_{02}b_{04}^d)$.

By Fact~\factassocses,
$\Ass(R/L) \subseteq \Ass(R/(L : c_{11})) \cup \Ass(R/(L + (c_{11})))$.
It will be proved that the only embedded prime of $J$
in this larger union set
$\Ass(R/(L : c_{11})) \cup \Ass(R/(L + (c_{11})))$
are the $Q_{4r\alpha\beta\gamma}$ or the $Q_{4,2\alpha\beta}$.

First of all,
$$
L + (c_{11}) = L_0 + \C_1 +
(s-fb_{01}^d,b_{01}^d-b_{02}^d, b_{04}^d-b_{03}^d, b_{01}-b_{04})
= L_0 + p_1,
$$
which equals the intersection of minimal components $p_{1\alpha\beta}$
if $L_0 = (c_{01}-c_{02}b_{01}^d, c_{03}-c_{02}, c_{04}-c_{02}b_{04}^d)$,
and is not associated to $J$ by Theorem~{\thmnox}\ otherwise.

Thus it remains to find the associated primes of $L : c_{11}$
in order to find the associated primes of $L$
which are also associated to $J$.
For this first note that $\JJ'' = \D_1 + c_{11}\JJ'''$
for some (obvious) ideal $\JJ'''$ in $R$.
Thus
$$
\eqalignno{
L :c_{11} &= L_0 + \D_1 + \JJ''' \cr
&\hskip2em+
(s-fb_{01}^d,b_{01}^d-b_{02}^d, b_{04}^d-b_{03}^d, b_{01}-b_{04},
b_{02}-b_{1i}b_{03}, c_{11}(b_{1i}-b_{1j}),1-b_{1i}^d). \cr
}
$$
Note that
$L : c_{11}b_{03}$ equals
$$
L_0 + \D_1 + \JJ''' +
(s-fb_{01}^d,b_{01}^d-b_{02}^d, b_{04}^d-b_{03}^d, b_{01}-b_{04},
b_{02}-b_{1i}b_{03}, b_{1i}-b_{1j},1-b_{1i}^d),
$$
which decomposes:
$$
\eqalignno{
&= \bigcap_{r = 2}^n \Bigl(L_0 + \D_1 + \cdots + \D_{r-1} + \C_r+ B_{r-1} \cr
&\hskip2em +
(s-fb_{01}^d,b_{01}^d-b_{02}^d,b_{04}^d-b_{03}^d, b_{01}-b_{04},
b_{02}-b_{1i}b_{03}, b_{1i}-b_{1j},1-b_{1i}^d)\Bigr) \cr
&= \bigcap_{r = 2}^n \Bigl(L_0 + p_r\Bigr). \cr
}
$$
As before,
when $L_0 = (c_{01}-c_{02}b_{01}^d, c_{03}-c_{02}, c_{04}-c_{02}b_{04}^d)$,
the above is just the intersection of some minimal components of $J$,
and when $L_0 = \C_0$,
the associated primes are of the form $\C_0 + P_{r\alpha\beta}$,
$r \ge 2$,
whence are not associated to $J$ by Theorems~\thmnox\ and~{\thmthirdemb}.

Thus it remains to find the associated primes of
$$
(L : c_{11}) + (b_{03})
= L_0 + \D_1 + \JJ''' +
(s,b_{01}^d, b_{01}-b_{04},b_{02},b_{03}, c_{11}(b_{1i}-b_{1j}),1-b_{1i}^d),
$$
which similarly decomposes (first add $c_{11}$ and colon with $c_{11}$) as
$$
\eqalignno{
&= \bigcap_{r = 2}^n
\Bigl(L_0 +\C_1 + D_2 + \cdots D_{r-1} + \C_r + B_{3,r-1}
+ (s,b_{01}^d, b_{01}-b_{04},b_{02}, b_{03},1-b_{1i}^d) \cr
&\hskip2em
+(b_{11}-b_{14})
+ (b_{12}-b_{2i}b_{13}, b_{2i}-b_{21}) \delta_{r > 2}
+ (b_{12}-b_{13}) \delta_{n = 2} \Bigr) \cr
&\hskip1em
\bigcap_{r = 2}^n \Bigl(L_0 + \D_1 + \cdots + \D_{r-1} + \C_r+ B_{r-1}
+ (s,b_{01}^d, b_{01}-b_{04},b_{02}, b_{03},b_{1i}-b_{1j},
1-b_{1i}^d) \Bigr), \cr
}
$$
from which it is easy to read off the associated primes.
By Theorems~\thmnox\ and~{\thmthirdemb},
only the $Q_{4r\alpha\beta\gamma}$ or the $Q_{4,2\alpha\beta}$
among these are embedded primes of $J$.

This proves the following:

\thm
\label{\thmoneleft}
The set of embedded primes of $J$ is contained in
$$
\eqalignno{
&\lbrace Q_{1\Lambda}, Q_{2\Lambda'\alpha}, Q_{3\Lambda'} |
\Lambda, \Lambda' \subseteq \{1,2,3,4\}, |\Lambda'| > 0 \rbrace \cr
&\hskip2em
\cup \lbrace
Q_{4,2\alpha\beta}\delta_{n=2},
Q_{4r\alpha\beta\gamma} \delta_{n>2} |
r = 2, \ldots, n; \alpha^d =\beta^d =\gamma^d = 1,
|\{\alpha, \beta, \gamma\}| > 1 \rbrace, \cr
&\hskip2em
\cup \Ass \left( {R \over (J : sc_{02})\hskip-0.1em+\hskip-0.1em(c_{02},f)}
\right), \cr
}
$$
where the explicitly listed $31 + 15d + d^2 \delta_{n = 2} + (n-1)(d^3-d)\delta_{n>2}$
prime ideals are indeed associated to $J = J(n,d)$.
\qed
\endb

Note that $(J : sc_{02}) + (c_{02},f)$
equals $K(n,d) + \C_0 + (s,f)$,
for an ideal $K(n,d)$ whose generators (given explicitly below)
do not involve any of the variables $s,f,c_{01}, c_{02}, c_{03}, c_{04}$:
$$
\eqalignno{
&g_{01} = b_{01} b_{03}^d - b_{04} b_{02}^d, \cr
& g_{1i} = c_{1i} \left(b_{02}-b_{1i} b_{03}\right), i = 1, \ldots, 4, \cr
& g_{1,4+i} = c_{1i} \left(b_{01}-b_{1i}^d b_{04}\right), i = 1, \ldots, 4, \cr
& g_{1ij} = c_{1i} c_{1j} \left(b_{1i} - b_{1j}\right), 1 \le i < j \le 4, \cr
& g_{21} = b_{04}^d c_{11} - b_{01}^d c_{12}, \cr
& g_{22} = b_{04}^d c_{14} - b_{01}^d c_{13},\cr
& g_{23} = b_{01}^d (c_{12} -c_{13}), \cr
& g_{24} = b_{04}^d (c_{12} b_{11} -c_{13} b_{14}), \cr
& g_{2,4+i} = b_{04}^d c_{12} c_{2i}
\left( b_{12}-b_{2i} b_{13} \right),
i = 1, \ldots, 4, \hbox{\ when $n > 2$}, \cr
& g_{25} = b_{04}^d c_{12} c_{2i}\left( b_{12}-b_{13} \right),
\hbox{\ when $n = 2$}, \cr
&g_{r1} = b_{01}^d c_{11} \cdots c_{r-3,1} \left(
c_{r-2,4} c_{r-1,1} - c_{r-2,1} c_{r-1,2}
\right), r = 2, \ldots, n, \cr
& g_{r2} = b_{01}^d c_{11} \cdots c_{r-3,1} \left(
c_{r-2,4} c_{r-1,4} - c_{r-2,1} c_{r-1,3}
\right), r = 2, \ldots, n, \cr
& g_{r3} = b_{01}^d c_{11} \cdots c_{r-2,1}
\left( c_{r-1,3}-c_{r-1,2} \right), r = 2, \ldots, n, \cr
& g_{r4} = b_{01}^d c_{11} \cdots c_{r-3,1} c_{r-2,4}
\left( c_{r-1,2} b_{r-1,1}-c_{r-1,3} b_{r-1,4} \right), r = 2, \ldots, n, \cr
&g_{r,4+i} = b_{01}^d c_{11} \cdots c_{r-3,1} c_{r-2,4} c_{r-1,2} c_{ri}
\left( b_{r-1,2}-b_{ri} b_{r-1,3} \right),
i = 1, \ldots, 4, r = 2, \ldots, n-1, \cr
&g_{n5} = b_{01}^d c_{11} \cdots c_{n-3,1} c_{n-2,4} c_{n-1,2}
\left( b_{n-1,2}-b_{n-1,3} \right). \cr
}
$$

The family of ideals $K(n,d)$ is analyzed in [S3].
In particular,
it is proved in [S3]
that this family also satisfies the doubly exponential ideal membership property.
Furthermore,
the set of associated primes of $K(n,d)$ recursively depends
on the set of associated primes of $K(n-1,d^2)$.

By Fact~\factfewvars,
any prime ideal associated to $K(n,d)$,
after adding $\C_0 + (s,f)$,
is possibly associated to $J(n,d)$.
In [S3] the obtained set of prime ideals possibly associated to $K(n,d)$
consists of 20 variously subscripted families
and the ideals associated to $K(n-1,d^2) + \C_1 + (b_{01}, b_{02}, b_{03}, b_{04})$,
where $K(n-1,d^2)$ involves the variables $c_{ri}, r \ge 2$,
and $b_{ri}, r \ge 1$.
From these families by Fact~\factfewvars\ %
then one easily constructs the corresponding families of prime ideals,
here subscripted with $5$ through $24$,
which are possibly associated to $J(n,d)$.
To list these families,
as usual,
$\Lambda$ always varies over all the subsets of $\{1,2,3,4\}$
and
$\Lambda'$ varies over all the non-empty subsets of $\{1,2,3,4\}$.
Also,
we will use the ideals
$$
T_r = (s,f) + \C_0 + \cdots + \C_r + (b_{ti} | i = 1, \ldots, 4; t = 0, \ldots, r-1).
$$
With this then the list of prime ideals in [S3] of prime ideals possibly associated
to $K(n,d)$ lifts to the following prime ideals
possibly associated to $J(n,d)$:
$$
\eqalignno{
&Q_{5r\Lambda'} = T_r +
\left(b_{r1},b_{r4}\right)
+ (c_{r+1,i} | i \not \in \Lambda) + \left(b_{r2}-b_{r+1,i}b_{r3},
b_{r+1,i}-b_{r+1,j} | i, j \in \Lambda \right), \cr
 &\hskip 3em \hbox{height } 8r+12, 0 \le r \le n-2, \cr
&Q_{6r} =T_r +\C_{r+1} + (b_{r1}b_{r3}^{d^{2^r}}-b_{r4}b_{r2}^{d^{2^r}}),
\hbox{height } 8r+11, 0 \le r \le n-2, \cr
&Q_{7r} = T_r +(c_{r+1,1}, c_{r+1,2},c_{r+1,4},b_{r1},b_{r2},b_{r+1,3},b_{r+1,4}), \cr
 &\hskip 3em \hbox{height } 8r+13, 0 \le r \le n-2, \cr
&Q_{8r\Lambda}=T_r +
(c_{r+1,1}, c_{r+1,4}, b_{r1}, b_{r2}, b_{r+1,2},b_{r+1,3},
c_{r+1,2}b_{r+1,1}-c_{r+1,3}b_{r+1,4}) \cr
&\hskip3em
+ (c_{r+2,i} | i \not \in \Lambda) + (1-b_{r+2,i} | i \in \Lambda),
\hbox{height } 8r+17, 0 \le r \le n-3, \cr
&Q_{9r}=T_r +
(c_{r+1,1}, c_{r+1,4}, b_{r1}, b_{r2}, b_{r+1,2}, b_{r+1,3},
c_{r+1,2}b_{r+1,1}-c_{r+1,3}b_{r+1,4}), \cr
 &\hskip 3em \hbox{height } 8r+13, 0 \le r \le n-2, \cr
&Q_{10r\Lambda}=T_r +
(c_{r+1,1}, c_{r+1,3}, c_{r+1,4}, b_{r1}, b_{r2}, b_{r+1,1}, b_{r+1,2}) \cr
&\hskip3em
+ (c_{r+2,i} | i \not \in \Lambda) + (b_{r+2,i} | i \in \Lambda),
\hbox{height } 8r+17, 0 \le r \le n-3, \cr
&Q_{11r\Lambda'}=T_r +
(c_{r+1,1}, c_{r+1,3}, c_{r+1,4}, b_{r1}, b_{r2}, b_{r+1,1}, b_{r+1,2},b_{r+1,3}) \cr
&\hskip3em
+(c_{r+2,i} | i \not \in \Lambda') +(b_{r+2,i}-b_{r+2,j} | i,j \in \Lambda'),
\hbox{height } 8r+17, 0 \le r \le n-3, \cr
&Q_{12r\Lambda\alpha}=T_r +
\C_1 + (b_{01}, b_{02}, b_{03}, b_{12},b_{13})
+ (c_{2i} | i \not \in \Lambda) + (b_{2i}-\alpha | i \in \Lambda),
\alpha^{d^{2^r}} = 1, \cr
 &\hskip 3em \hbox{height } 8r+19, 0 \le r \le n-3, \cr
&Q_{13r\Lambda'}=T_r +
\C_{r+1} + (b_{r1}, b_{r2}, b_{r3}, b_{r+1,2},b_{r+1,3}) \cr
&\hskip3em
+ (c_{r+2,i} | i \not \in \Lambda')
+ (b_{r+2,i}-b_{r+2,j} | i,j \in \Lambda'),
\hbox{height } 8r+18, 0 \le r \le n-3, \cr
&Q_{14r\Lambda}=T_r +\C_{r+1} +
(b_{r1}, b_{r2}, b_{r3}, b_{r4}, b_{r+1,2}, b_{r+1,3})
+ (c_{r+2,i} | i \not \in \Lambda)  \cr
&\hskip3em+
(b_{r+2,i}-b_{r+2,j} | i,j \in \Lambda),
\hbox{height } 8r+19 + \delta_{\lambda = \emptyset}, 0 \le r \le n-3, \cr
&Q_{15r\Lambda} = T_r +(c_{r+1,1},c_{r+1,3}-c_{r+1,2}, c_{r+1,4},
b_{r1},b_{r2},b_{r+1,1},b_{r+1,2},b_{r+1,3},b_{r+1,4}) \cr
&\hskip3em
+ (c_{r+2,i} | i \not \in \Lambda)
+ (1-b_{r+2,i} | i \in \Lambda),
\hbox{height } 8r+19, 0 \le r \le n-3, \cr
&Q_{16r\Lambda} = T_r +\C_{r+1} +
(b_{r1},b_{r2},b_{r+1,1},b_{r+1,2},b_{r+1,3},b_{r+1,4}) \cr
&\hskip3em
+ (c_{r+2,i} | i \not \in \Lambda)+ (1-b_{r+2,i} | i\in \Lambda),
\hbox{height } 8r+20, 0 \le r \le n-3, \cr
&Q_{17r\Lambda} = T_r +\C_{r+1} +
(b_{r1},b_{r2},b_{r+1,1},b_{r+1,2},b_{r+1,3},b_{r+1,4})
+ (c_{r+2,i} | i \not \in \Lambda) \cr
&\hskip3em
+ (b_{r+2,i}-b_{r+2,j} | i, j \in \Lambda),
\hbox{height } 8r+19 + \delta_{\lambda = \emptyset}, 0 \le r \le n-3, \cr
&Q_{18r\Lambda} = T_r +\C_{r+1} +
(b_{r1},b_{r2},b_{r3},b_{r+1,1},b_{r+1,2},b_{r+1,3},b_{r+1,4})
+ (c_{r+2,i} | i \not \in \Lambda) \cr
&\hskip3em
+ (b_{r+2,i}- b_{r+2,j} | i,j \in \Lambda),
\hbox{height } 8r+20 + \delta_{\lambda = \emptyset}, 0 \le r \le n-3, \cr
&Q_{19r\Lambda'\alpha} = T_r +\C_{r+1} +
(b_{r1},b_{r2},b_{r3},b_{r+1,1},b_{r+1,2},b_{r+1,3},b_{r+1,4}) \cr
&\hskip3em
+ (c_{r+2,i} | i \not \in \Lambda')
+ (b_{r+2,i}- \alpha | i \in \Lambda'), \alpha^{d^{2^{r+1}}} = 1,
\alpha^{d^{2^r}} \not= 1, \cr
 &\hskip 3em \hbox{height } 8r+21, 0 \le r \le n-3, \cr
&Q_{20\Lambda'\alpha} = T_r +\C_{r+1} +
(b_{r1},b_{r2},b_{r3},b_{r+1,1},b_{r+1,2},b_{r+1,3},b_{r+1,4}) \cr
&\hskip3em
+ (c_{r+2,i} | i \not \in \Lambda')
+ (b_{r+2,i}- \alpha | i \in \Lambda'), \alpha^{d^{2^r}} = 1, \cr
 &\hskip 3em \hbox{height } 8r+21, 0 \le r \le n-3, \cr
&Q_{21rt}= T_r +\D_{r+2} + \cdots + \D_{t-1} + \C_t
+ B_{2t-1} \cr
&\hskip3em+ (c_{r+1,1}-b_{r+1,2}^{d^2}c_{r+1,2},
c_{r+1,4}-c_{r+1,1}, c_{r+1,3}-c_{r+1,2}) \cr
&\hskip3em
+(b_{r2}-b_{r+1,2}b_{r3},
b_{r+1,2}-b_{r+1,i},b_{r1}-b_{r+1,2}^db_{r4} ), \cr
 &\hskip 3em \hbox{height } 7t + r + 4 \delta_{t < n}, 0 \le r \le n-2, \cr
&Q_{22rt}= T_r +\C_{r+1}+ D_{r+2} + \cdots + \D_{t-1} + \C_t+ B_{2t-1} \cr
&\hskip3em
+ \left( b_{r2}-b_{r+1,2}b_{r3},
b_{r+1,2}-b_{r+1,i},
b_{r1}-b_{r+1,2}^db_{r4} \right), \cr
 &\hskip 3em \hbox{height } 7t + r + 1 + 4 \delta_{t < n}, 0 \le r \le n-2, \cr
&Q_{23r,n-2\alpha\beta} = T_r +\C_{r+1}+ \C_{r+2}
+ \left(b_{r1}-b_{r+1,2}^{d^{2^r}} b_{r4},b_{r2}, b_{r3}, b_{r+1,1}-b_{r+1,4}
\right)\cr
&\hskip3em
+ \left(b_{r+1,2}-\alpha b_{r+1,3},b_{r+1,1}-\beta b_{r+1,3}\right),
\alpha^{d^{2^r}} = \beta^{d^{2^r}} = 1,
\hbox{height } 8n, \cr
&Q_{23rt\alpha\beta} = T_r +\C_{r+1}+ \D_{r+2} + \cdots + \D_{t-1}
+ \C_t + B_{3t-1}
+ \left(b_{r1}-b_{r+1,2}^db_{r4},b_{r2}, b_{r3}\right) \cr
&\hskip 3em
+ \left( b_{r+1,2}-\alpha b_{r+1,3},b_{r+1,1}-\beta b_{r+1,3},
 b_{r+1,1}-b_{r+1,4},
 b_{r+2,i}-\alpha
\right), \cr
&\hskip 3em \alpha^{d^{2^r}} = \beta^{d^{2^r}} = 1, \hbox{height } 7t + r + 2 + 4 \delta_{t < n}, 0 \le r \le n-3, \cr
&Q_{24} = T_{n-1} + (b_{n-1,1} - b_{n-1,4}, b_{n-1,2} - b_{n-1,3}),
 \hbox{height } 8n. \cr
}
$$
\vskip2em

Thus finally:

\thm
\label{\thmfinal}
With $n \ge 2$,
the set of embedded primes of the Mayr-Meyer ideal $J = J(n,d)$ is contained in
the set
$$
\eqalignno{
&\lbrace Q_{1\Lambda}, Q_{2\Lambda'\alpha},Q_{3\Lambda'}, Q_{24} \rbrace
\bigcup_{r = 2}^n
\lbrace
 Q_{4,2\alpha\beta} \delta_{n = 2},
 Q_{4r\alpha\beta\gamma} \delta_{n > 2}
| \alpha^d = \beta^d = \gamma^d = 1,
|\{\alpha, \beta, \gamma\}| > 0 \rbrace \cr
&\hskip2em\bigcup_{r = 0}^{n-2}
\{Q_{5r\Lambda'}, Q_{jr}, Q_{krt} |j = 6,7,9; k = 21,22;t = r+2, \ldots, n\} \cr
&\hskip2em
\bigcup_{r = 0}^{n-3}
\{Q_{ir\Lambda}, Q_{kr\Lambda'}|
i = 8, 10, 14,15,16,17,18; k = 11, 13\} \cr
&\hskip2em
\bigcup_{r = 0}^{n-3}
\{Q_{12r\Lambda\alpha},Q_{19r\Lambda'\alpha'},Q_{20r\Lambda'\alpha}
| \alpha^{d^{2^r}}=1, \alpha'^{d^{2^{r+1}}}= 1, \alpha'^{d^{2^r}} \not = 1\}\cr
&\hskip2em \bigcup_{r = 0}^{n-3}
\lbrace Q_{23rt\alpha\beta} | t = r+2, \ldots, n;
\alpha^{d^{2^r}}=\beta^{d^{2^r}}=1\rbrace
\bigcup
\lbrace Q_{23,n-2,n,1\alpha} |
\alpha^{d^{2^{n-2}}}=1\rbrace. \cr
}
$$
where $\Lambda$ varies over all the subsets of $\{1,2,3,4\}$,
and $\Lambda'$ varies over all the non-empty subsets of $\{1,2,3,4\}$.
\qed
\endb

\remark
It was proved in Section~\sectmoreemb\ that
the $Q_{1\Lambda}, Q_{2\Lambda'\alpha},Q_{3\Lambda'}$
are indeed associated to $J$,
and in Section~\sectmoremoreemb\ %
the same was proved for the $Q_{4r\alpha\beta\gamma}$
and the $Q_{42\alpha\beta}$.
\endb

The last theorem proves
that the Mayr-Meyer ideal $J(n,d)$ has at most \embprno\ embedded prime ideals.

Also,
for all $n \ge 2$,
none of the maximal ideals is associated to the Mayr-Meyer ideals.

Whereas the theorem above gives some information on the structure
of the associated prime ideals of $J(n,d)$,
much is left to be done to answer the Bayer-Huneke-Stillman question.
I end this paper with a list of questions:
\item{1.}
Some of the prime ideals in Theorem~\thmfinal\ may not be associated to $J(n,d)$.
Find all such primes,
or in other words,
find the exact set of embedded primes of $J(n,d)$,
not just a set containing it.
In particular,
determine if the set of associated primes of $J(n,d)$
is truly doubly exponential in $n$.
\item{2.}
Determine if any of the associated prime ideals of $J(n,d)$
play a crucial role in the doubly exponential behavior.
The prime ideals $Q_{23,n-2,n,1,\alpha}$ and $Q_{24}$ may be likely candidates.
\item{3.}
The ideal $J(n,d) + (s,f)^2 +
\sum_{r=0}^{n-1} (c_{r1},c_{r2},c_{r3},c_{r4})^2$
exhibits the same doubly exponential syzygetic behavior as $J(n,d)$.
It has height $2 + 4n$,
whereas $J(n,d)$ has height~2.
What kind of primary decomposition or associated prime ideal
structure does this larger ideal exhibit?

\vskip 4ex

\bigskip
\leftline{\bf References}
\bigskip

\bgroup
\font\eightrm=cmr8 \def\rm{\fam0\eightrm}
\font\eightit=cmti8 \def\it{\fam\itfam\eightit}
\font\eightbf=cmbx8 \def\bf{\fam\bffam\eightbf}
\font\eighttt=cmtt8 \def\tt{\fam\ttfam\eighttt}
\rm
\baselineskip=9.9pt
\parindent=3.6em

\item{[BS]}
D.\ Bayer and M.\ Stillman,
On the complexity of computing syzygies,
{\it J.\ Symbolic Comput.}, {\bf 6} (1988), 135-147.

\item{[D]}
M.\ Demazure,
Le th\'eor\`eme de complexit\'e de Mayr et Meyer,
{\it G\'eom\'etrie alg\'ebrique et applications, I
(La R\'abida, 1984)}, 35-58,
Travaux en Cours, 22, {\it Hermann, Paris}, 1987.

\item{[GS]}
D.\ Grayson and M.\ Stillman,
Macaulay2. 1996.
A system for computation in algebraic geometry and commutative algebra,
available via anonymous {\tt ftp} from {\tt math.uiuc.edu}.

\item{[GPS]}
G.-M.\ Greuel, G.\ Pfister and H.\ Sch\"onemann,
Singular. 1995.
A system for computation in algebraic geometry and singularity theory.
Available via anonymous {\tt ftp} from {\tt helios.mathematik.uni-kl.de}.

\item{[H]}
G.\ Herrmann,
Die Frage der endlich vielen Schritte in der Theorie der Polynomideale,
{\it Math.\ Ann.}, {\bf 95} (1926), 736-788.

\item{[K]}
J.\ Koh,
Ideals generated by quadrics exhibiting double exponential degrees,
{\it J.\ Algebra}, {\bf 200} (1998), 225-245.

\item{[MM]}
E.\ Mayr and A.\ Meyer,
The complexity of the word problems for commutative semigroups
and polynomial ideals,
{\it Adv.\ Math.}, {\bf 46} (1982), 305-329.

\item{[S1]}
I.\ Swanson,
The first Mayr-Meyer ideal,
preprint, 2001.

\item{[S2]}
I.\ Swanson,
The minimal components of the Mayr-Meyer ideals,
preprint, 2002.

\item{[S3]}
I.\ Swanson,
A new family of ideals with the doubly exponential ideal membership property,
preprint, 2002.

\egroup

\vskip 4ex
\noindent
{\sl
New Mexico State University - Department of Mathematical Sciences,
Las Cruces, New Mexico 88003-8001, USA.
E-mail: {\tt iswanson@nmsu.edu}.
}

\end